\documentclass[11pt]{article}
\usepackage{epsfig,numbysec}
\usepackage{amsmath}
\usepackage{amssymb}
\hsize \textwidth
\advance \hsize by -\marginparwidth
\evensidemargin \oddsidemargin
\usepackage{amssymb}
\advance\hoffset by 5mm

\newcommand{\nit}{\noindent}
\newcommand{\no}{\nonumber}
\newcommand{\be}{\begin{equation}}
\newcommand{\ee}{\end{equation}}
\newcommand{\bi}{\begin{itemize}}
\newcommand{\ei}{\end{itemize}}
\newcommand{\br}{\begin{eqnarray}}
\newcommand{\er}{\end{eqnarray}}

\newcommand{\eps}{\varepsilon}

\newcommand{\qed}{$\Box$}

\newcommand{\defeq}{\stackrel{\Delta}{=}}

\newcommand{\abs}[1]{\lvert #1 \rvert}
\newcommand{\norm}[1]{\lVert #1 \rVert}

\newcommand{\commentout}[1]{}

\newtheorem{thm}[theorem]{Theorem}
\newtheorem{prop}[theorem]{Proposition}
\newtheorem{lemma}[theorem]{Lemma}
\newtheorem{rem}[theorem]{Remark}
\newtheorem{definition}[theorem]{Definition}
\newtheorem{corollary}[theorem]{Corollary}
\newcommand{\Rm}{{\mathbb R}}
\newcommand{\Nm}{{\mathbb N}}
\newcommand{\Em}{{\mathbb E}}
\newcommand{\Pm}{{\mathbb P}}

\newcommand{\pdr}[2]{\frac{\partial{#1}}{\partial{#2}}}

\begin{document}
\numberbysection
\title{Traveling waves in a one-dimensional random medium}

\author{James Nolen\thanks{Department of Mathematics, Stanford University,
Stanford, CA 94305, USA. (nolen@math.stanford.edu).} \;
and \; Lenya Ryzhik
\thanks{Department of Mathematics, University of Chicago, Chicago,
IL 60637, USA.
(ryzhik@math.uchicago.edu).}}

\date{\today}

\maketitle

\begin{abstract}
We consider solutions of a scalar reaction-diffusion equation of the ignition type
with a random, stationary and ergodic reaction rate. We show that solutions of the
Cauchy problem spread with a deterministic rate in the long time limit. We also
establish existence of generalized random traveling waves and of transition fronts
in general heterogeneous media.
\end{abstract}

\section{Introduction}

We consider solutions to the equation
\be
u_t = \Delta u + f(x,u,\omega), \;\; x \in \mathbb{R} \label{eq0}
\ee
where $f(x,u,\omega)$ is a random ignition-type nonlinearity that is
stationary with respect to translation in $x$. The function $f$ has
the form $f(x,u,\omega) = g(x,\omega)f_0(u)$. Here, $f_0(u)$ is an
ignition-type nonlinearity with an ignition temperature $\theta_0 \in (0,1)$:
$f_0(u)$ is a Lipschitz function,  and, in addition,
\[
\hbox{$f_0(u) = 0$ for $u \notin (\theta_0,1)$,
$f_0(u) > 0$ for $u \in (\theta_0,1)$, $f'(1) < 0$.}
\]
The reaction rate $g(x,\omega)$, $x \in \mathbb{R}$, is a stationary,
ergodic random field defined over a probability space $(\Omega,\Pm,
\mathcal{F})$:   there exists a group $\{ \pi_x \}$, $x\in\Rm$, of
measure-preserving transformations acting ergodically on $(\Omega,
\Pm, \mathcal{F})$ such that $g(x + h,\omega) = g(x,\pi_h \omega)$.
We suppose that $g(x,\omega)$ is almost surely Lipschitz continuous with
respect to $x$ and that there are deterministic constants $g^{min},
g^{max}$ such that
\[
0 < g^{min} \leq g(x,\omega) \leq g^{max} < \infty
\]
holds almost surely. Thus, we have
\[
f^{min}(u) \leq f(x,u,\omega) \leq
f^{max}(u),
\]
where $f^{min}(u) = g^{min}f_0(u)$ and $f^{max}(u) = g^{max}f_0(u)$
are both ignition-type nonlinearities with the same ignition
temperature. We assume that the probability space $\Omega =
C(\mathbb{R};[g^{min},g^{max}])$ and that $\mathcal{F}$ contains the
Borel $\sigma$-algebra generated by the compact open topology (the
topology of locally uniform convergence) on
$C(\mathbb{R};[g^{min},g^{max}])$.

We are interested in the following two issues: first, how do solutions of the
Cachy problem for (\ref{eq0}) with a compactly supported non-negative initial data
spread in the long time limit? Second, do there exist special solutions
of (\ref{eq0}) that generalize the notion of a traveling front in the homogeneous case?

It is well known since the pioneering work
by Ya. Kanel~\cite{Kanel} that in the uniform case:
\begin{equation}\label{intro-homogen}
u_t=\Delta u+f(u)
\end{equation}
with an ignition-type nonlinearity $f(u)$, all solutions with the initial data
$u_0(x)=u(0,x)$ in a class $I \subset C_c(\Rm)$,
$0 \leq u_0(x)\leq 1$, propagate
with the same speed $c^*$ in the sense that
\begin{equation}\label{intro-spread-homogen0}
\lim_{t\to+\infty} u(t,ct)=0\hbox{ for $|c|>c^*$},
\end{equation}
and
\begin{equation}\label{intro-spread-homogen1}
\lim_{t\to+\infty} u(t,ct)=1\hbox{ for $|c|<c^*$}.
\end{equation}
The initial data is restricted to the class $I$ to preclude the possibility of the so-called quenching phenomenon where $u \to 0$ uniformly in $x$ as $t \to \infty$. In particular,
$I$ contains functions that are larger than
$\theta_0 + \epsilon$ on a sufficiently large interval, depending on $\epsilon > 0$.
The constant $c^*$ above is the speed of the unique traveling wave solution
$u(t,x)=U(x-c^*t)$ of (\ref{intro-homogen}):
\[
-c^*U'=U''+f(U),~~U(-\infty)=1,~~U(+\infty)=0.
\]

As far as heterogeneous media are concerned this result has been
extended to the periodic case: J.~Xin \cite{Xin-ARMA,Xin-JSP}, and
H.~Berestycki and F.~Hamel \cite{BH} have shown that when the function
$f(x,u)$ is periodic in $x$, equation (\ref{eq0}) admits special solutions of
the form $u(t,x)=U(x-c^*t,x)$, called pulsating fronts, which are
periodic in the second variable and satisfy
\[
\hbox{$U(s,x)\to 1$ as $s\to-\infty$, and $U(s,x)\to 0$ as $x\to +\infty$.}
\]
H.~Weinberger~\cite{Weinberger} has proved that solutions with general
non-negative compactly supported initial data spread with the speed
$c^*$ in the sense of (\ref{intro-spread-homogen0})-(\ref{intro-spread-homogen1}),
though the
spreading rates to the left and right may now be different.

The purpose of the present paper is to extend the result
of~\cite{Weinberger} to the stationary random ergodic case, and show
that special solutions which generalize the notion of a pulsating
front to random media exist.

\subsection*{Deterministic spreading rates}

Our first result concerns the asymptotic behavior of
solutions to the Cauchy problem for (\ref{eq0}) with compactly
supported initial data. We show that for sufficiently large initial
data the solution develops two diverging fronts that propagate with a
deterministic asymptotic speed. Specifically, we prove the following.
\begin{thm} \label{theo:asympspeed}
Let $w(x,t,\omega)$ solve (\ref{eq0}) with compactly supported deterministic initial
data $w_0(x)$, $0\le w_0(x)\le 1$. Let $h \in (\theta_0,1)$ and suppose that
$w_0 \geq h$ on an interval of size $L > 0$. There exist deterministic constants
$c^*_- < 0 < c^*_+$ such that for any $\eps > 0$, the limits
\[
\lim_{t \to \infty} \inf_{c \in [c^*_- + \eps,c^*_+ - \eps]}
w(t,ct,\omega) = 1
\]
and
\[
\lim_{t \to \infty} \sup_{c \in (-\infty,c^*_- - \epsilon]
\cup [c^*_+ + \epsilon,\infty)}  w(t,ct,\omega) = 0
\]
hold almost surely with respect to $\Pm$, if $L$ is sufficiently
large. The constants $c^*_- , c^*_+$ are independent of $h$ and $L$.
\end{thm}
The condition that $L$ be sufficiently large is necessary only to
exclude the possibility of uniform convergence to zero \cite{Kanel}.

Using large deviation techniques, Freidlin, and Freidlin and G\"artner
(see \cite{FR1}, section 7.4, \cite{FG,freidlin-3,freidlin-2,Freidlin-02}) proved a
similar asymptotic result in the case that $f_0(u)$ is of KPP-type
satisfying $f(u) \leq f'(0) u$ (e.g. $f_0(u) = u(1-u)$).  Moreover, the asymptotic
speed can be identified by a variational principle that arises from the linearized
problem at $u=0$. This asymptotic spreading result has been extended recently to time-dependent
random media in \cite{NX2,NX3}. The problem with a KPP nonlinearity also admits homogenization, both in the periodic \cite{MS} and
random \cite{KRV,LS} cases.
However, in all aforementioned papers,
the KPP condition $f(u) \leq f'(0) u$ seems to be essential, and
the techniques do not extend to the present case where $f$ vanishes
when $u$ is close to zero. To the best of our knowledge
Theorem~\ref{theo:asympspeed} is the first result on the deterministic spreading rates
of solutions of reaction-diffusion equations with a non-KPP nonlinearity in
a random medium.

\subsection*{Random traveling waves}

Two generalizations of the notion of a traveling front in a uniform medium
for general (non-periodic)
inhomogeneous media were proposed.
Shen in \cite{Shen1}, and Berestycki and Hamel in \cite{BH-CRAS,BH-Brezis}
have introduced generalized transition fronts
(called wave-like solutions in \cite{Shen1}) -- these are global in time solutions
that, roughly speaking, have an interface which ``stays together" uniformly in time.
On the other hand, H. Matano has defined a generalized traveling wave as a
global in time solution whose shape is ``a continuous function of the current
environment" \cite{Matano}. These notions are not equivalent: there exist
transition fronts of the KPP equation with constant coefficients that are
not traveling waves in the usual sense (and hence not generalized
traveling waves in the sense of Matano as there is only one environment in the case of
a uniform medium and thus only one solution profile) \cite{HN1,HN2}.

Matano's definition was formalized by W. Shen in \cite{Shen1} as follows.
\begin{definition}[see \cite{Shen1}, Def. 2.2] \label{def:rtw}
A solution $\tilde w(t,x,\omega):\mathbb{R} \times \mathbb{R} \times \Omega \to \mathbb{R}$
of (\ref{eq0}) is called a {\bf random traveling wave} if the following hold:
\bi
\item[(i)] For almost every $\omega \in \Omega$, $\tilde w(t,x,\omega)$ is a classical solution of
(\ref{eq0}) for all $t \in \mathbb{R}$.
\item[(ii)] The function $\tilde w(0,x,\omega)$ is measureable with respect to $\omega$.
\item[(ii)] $0 < \tilde w(0,x,\omega) < 1,\quad \forall x \in \mathbb{R}$.
\item[(iii)] $\lim_{x \to +\infty} \tilde w(0,x,\omega) = 0$.
\item[(iv)] $\lim_{x \to -\infty} \tilde w(0,x,\omega) = 1$.
\item[(v)]  There exists a measureable function $\tilde X(t,\omega):\mathbb{R} \times \Omega \to \mathbb{R}$ such that
\[
\tilde w(t,x,\omega) = \tilde w(0,x - \tilde X(t,\omega),\pi_{\tilde X(t,\omega)} \omega).
\]
\ei
The function $W(x,\omega):\mathbb{R} \times \Omega \to \mathbb{R}$ defined
by $W(x,\omega) = \tilde w(0,x,\omega)$ is said to {\bf generate} the random
traveling wave.
\end{definition}
The random function $W(x,\omega)$ is the profile of the wave in the moving
reference frame defined by the current front position $\tilde X(t,\omega)$.
In the pioneering paper \cite{Shen1}, Shen has established some general criteria  for the
existence of a traveling wave in ergodic spatially and temporally varying media
and also proved some important properties of the wave.  In particular, as
shown in \cite{Shen1} (see Theorem~B, therein),
Definition \ref{def:rtw} of a random traveling wave generalizes the notion of a pulsating
traveling front in the periodic case.
More precisely, if $f$ is actually periodic in $x$, then a
random traveling wave solution of (\ref{eq0}) is a pulsating traveling
front solution in the sense of \cite{BH,Xin-ARMA,Xin-JSP}.

However, the only example provided in \cite{Shen1} where the results of \cite{Shen1} ensure existence
of a random traveling wave is a bistable reaction-diffusion equation of the form
\[
u_t=u_{xx}+(1-u)(1+u)(u-a(t)),
\]
where $a(t)$ is a stationary ergodic random process. As far as we are aware, no
other examples of such traveling waves in non-periodic media have been exhibited.
In this paper we construct
a Matano-Shen traveling wave in a spatially varying
ergodic random medium for (\ref{eq0}) with
an ignition-type nonlinearity.
\begin{thm}\label{theo:rtwexist}
There exists a random traveling wave solution of (\ref{eq0}).
Moreover, $\tilde w(\tilde X(t,\omega),\omega) = \theta_0$ for all
$t \in \mathbb{R}$, $\tilde X'(t,\omega) \geq 0$, and
\be
\lim_{t \to \infty} \frac{\tilde X(t,\omega)}{t} = c^*_+ \label{rtwsamec}
\ee
holds almost surely, where $c^*_+$ is the same constant in Theorem \ref{theo:asympspeed}.
\end{thm}
Since $\tilde X(t,\omega)$ is increasing in $t$, we may define its inverse $\tilde T(x,\omega):\mathbb{R} \times \Omega \to \mathbb{R}$ by
\be
x = \tilde X(\tilde T(x,\omega),\omega).
\ee
This may be interpreted as the time at which the interface reaches the position $x \in \mathbb{R}$.  The following Corollary says that the statistics of the profile of the wave as the wave passes
through the point $\xi$ are invariant with respect to $\xi$:
\begin{corollary} \label{cor:profinv}
The function $\tilde w(x + \xi,\tilde T(\xi,\omega),\omega)$ is stationary with respect to shifts in $\xi$.
\end{corollary}
This a direct analog of the corresponding property of a pulsating front in a periodic
medium: the profile of a pulsating front at the time $T(\xi)$ it passes a point $\xi$
is periodic in $\xi$.

We believe
the present article gives the first construction of such a wave in a
spatially random medium. To construct the wave, we use a dynamic
approach from \cite{Shen1} combined with some analytical estimates
needed to show that the construction produces a nontrivial result.

\subsection*{Generalized transition fronts}

Our last result concerns existence of the transition fronts for (\ref{eq0})
in the sense of Berestycki and Hamel, and Shen,
in general heterogeneous (non-random)
media with the reaction rate  uniformly bounded from below and above.
Let us recall first the definition of a generalized transition wave.
\begin{definition}\label{def-tr-front}
A global in time solution $\tilde v(t,x)$, $t\in\Rm$, $x\in\Rm$, of (\ref{eq0})
is called a transition wave
if for any $h,k \in (0,1)$ with $h > k$, we have
\be
0 \leq \theta_{k}^+(t,\omega) - \theta_h^-(t,\omega) \leq C \label{wlsineq}
\ee
for all $t \in\Rm$, where
\br
\theta_h^-(t,\omega) = \sup \{ x \in \mathbb{R} |\; \tilde v(t,x',\omega) > h \;\;\forall x' < x \} \no \\
\theta_{k}^+(t,\omega)=  \inf \{ x \in \mathbb{R} |\; \tilde v(t,x',\omega) < k \;\;\forall x' > x \}
\er
and $C = C(h,k)$ is a constant independent of $t$ and $\omega$.
\end{definition}
Roughly speaking, a transition wave is a global in time solution for which
there are uniform, global-in-time bounds on the width of the interface.
Basic properties of transition waves were investigated in \cite{BH-CRAS,BH-Brezis}.

\begin{thm}\label{thm-transition}
Let $f(x,u)$ be a nonlinearity such that $g^{min}f_0(u)\le f(x,u)\le g^{max}f_0(u)$, with
the constants
$g^{min}>0$, $g^{max}<+\infty$ and $f_0(u)$ an ignition-type nonlinearity.
Then there exists a transition front solution $u(t,x)$, $t\in\Rm$, $x\in\Rm$,
of (\ref{eq0}) which
is monotonically increasing in time: $u_t(t,x)>0$.
In addition, there exists a unique point $X(t)$ such that $u(t,X(t))\equiv \theta_0$,
and a constant $p>0$ so that $u_x(t,X(t))<-p$ for all $t\in\Rm$.
\end{thm}

As this paper was written we learned about the concurrent work
by A.~Mellet and J.-M.~Roquejoffre \cite{MRoq}. They consider the free boundary
limit for (\ref{eq0}) in an ergodic random medium in the spirit of  \cite{CafMel}.
In particular, they also prove Theorem~\ref{thm-transition} as a necessary intermediate
step as well as other interesting results.

The paper is organized as follows. In Section \ref{sec:tran},
we study solutions to (\ref{eq0}) that are monotone increasing in
time and prove Theorem~\ref{thm-transition}. The main ingredients in the
proof are Propositions~\ref{thm:khseparation} and~\ref{thm:uxlower}
which show that the interface (the region where $\eps < u <1-\eps$, for some
$\eps > 0$) may not be arbitrarily wide and
must move forward with an instantaneous speed that is bounded above
and below away from zero. These estimates are also used later in the
proof of the asymptotic spreading and in the construction of the
random traveling waves.  In Section \ref{sec:asympspread} we prove
Theorem \ref{theo:asympspeed}, first for monotone increasing solutions
and then for general compactly supported data.  In Section
\ref{sec:rtw} we construct the random traveling wave and prove Theorem
\ref{theo:rtwexist} and Corollary \ref{cor:profinv}.

Throughout the paper we denote by $C$ and $K$ universal constants that
depend only on the constants $g^{min}$ and $g^{max}$, and the function
$f_0(u)$.

{\bf Acknowledgment.} Most of this work was done during the visit by
LR to Stanford University, LR thanks the Stanford Mathematics Department
for its hospitality. JN was supported by an NSF Postdoctoral fellowship, and
LR by the NSF grant DMS-0604687 and ONR. We thank
A. Mellet and J.-M. Roquejoffre for describing to us the results of \cite{MRoq}.

\section{Existence of a generalized transition front}\label{sec:tran}


\subsection*{Monotonic in time solutions}

In this section we prove Theorem~\ref{thm-transition}. The generalized
transition wave is constructed as follows. We consider a sequence of
solutions $u_n(t,x)$ of (\ref{eq0}) (with $f(x,u,\omega)$ replaced by $f(x,u)$ as in the statement of the theorem) defined for $t\ge -n$, with the
Cauchy data
\begin{equation}\label{trans-indata}
u_n(t=-n,x)=\zeta(x-x_0^n),~~\zeta(x):=\max(\hat\zeta(x),0).
\end{equation}
The choice of the initial shift $x_0^n$ is specified below while the function
$\hat\zeta(x)$ is positive on an open interval and is a
sub-solution for (\ref{eq0}):
\begin{equation}\label{trans-hatzeta}
-\hat \zeta''(x) = f^{min}(\hat \zeta(x))\le f(x,\hat\zeta),
\end{equation}
with $f^{min}(u)=g^{min}f_0(u)\le f(x,u)$.
It is constructed as follows.
For a given $h_0 \in (\theta_0,1)$ and $x \in \mathbb{R}$, let $\hat \zeta(x)$ satisfy
\[
-\hat \zeta''(x) = f^{min}(\hat \zeta(x)),~~\hat \zeta(0) = h_0,
~~\hat \zeta'(0) = 0,
\]
with the convention that $f^{min}(u) = 0$ for $u < 0$ above. To fix ideas we may
set $h_0=(1+\theta_0)/2$ in (\ref{trans-hatzeta}).
Let us define $z_1 = \min \{ x > 0 \; | \; \hat \zeta(x) = \theta_0 \}$ and
$z_2 = \min \{ x > 0 \; | \; \hat \zeta(x) = 0 \}$. The function $\hat \zeta$ satisfies
the following elementary properties
\bi
\item $\hat \zeta(-x) = \hat \zeta(x)$ for all $x \in \mathbb{R}$
\item $0 \leq \hat \zeta(x) \leq h_0 = \hat \zeta(0)$ for all $x \in [-z_2,z_2]$
\item $\hat \zeta(x)$ is strictly concave for $x \in (-z_1,z_1)$
\item $\hat \zeta(x)$ is linear for $x \in [-z_2,-z_1]$ and $x \in [z_1,z_2]$
\item $\hat \zeta(-z_2) = \hat\zeta(z_2) = 0$.
\ei

As in \cite{Mal-Roq,Rjfr} it follows that $u_n(t,x)\ge u_n(-n,x)$ for
$t\ge -n$, and $u_n(t,x)$ is monotonically increasing in time to $\bar
u\equiv 1$.
\begin{lemma} \label{lem:monzeta}
Let $u_n(t,x)$ solve (\ref{eq0}) with initial data
(\ref{trans-indata}) at time $t = -n$. Then, $u_n(t,x)$ is strictly
increasing in $t$:
\begin{equation}\label{trans-monotone}
\pdr{  u_n}{  t}(t,x) > 0\hbox{ for all $t > -n$},
\end{equation}
and, moreover,
\begin{equation}\label{trans-toone}
\lim_{t \to \infty} u_n(t,x) = 1\hbox{ locally uniformly in $x$.}
\end{equation}
\end{lemma}
{\bf Proof.} Since
\[
-\zeta_{xx} \leq f(x,\zeta),
\]
the maximum principle implies that $u_n(t,x) \ge u_n(-n,x) = \zeta(x)$
for all $t> -n$. Applying the maximum principle to the function
$w(t,x) = u_n(t + \tau,x) - u_n(t,x)$, for $\tau > 0$ fixed,
we see that, as $w(-n,x)\ge 0$, we have  $w(t,x) > 0$ for all
$t > -n$; thus, $u_n$ is monotonically increasing in time and (\ref{trans-monotone})
holds.

Since $u_n$ is monotone in $t$, the limit $\bar u(x) = \lim_{t \to \infty} u_n(t,x)$
exists and satisfies
\begin{equation}\label{trans-baru}
\hbox{$\bar u_{xx} = -f(x,\bar u)$, $0 < \bar u(x) \leq 1$,
$\max_x \bar u > h_0 > \theta_0$.}
\end{equation}
Note that $\bar u_{xx}  = 0$ on the set $\{ \bar u < \theta_0\}$, so $\bar u$ is linear
there. It is easy to see that this implies that this set must be empty because of
the lower bound  $\bar u \geq 0$ and the fact that $\max_x \bar u > h_0 > \theta_0$.
Hence,  we have $\bar u \geq \theta_0$.

Now, (\ref{trans-baru}) implies that $\bar u$ is concave. Since
$\theta_0 \leq \bar u \leq 1$, this implies $\bar u$ is constant, so
that $f(x,\bar u) = -\bar u_{xx} \equiv 0$. This fact and $\max_x \bar
u > h_0$ implies that $\bar u \equiv 1$. The local uniformity of the
limit follows from standard regularity estimates for $u$. $\Box$

\subsection*{The initial shift}

The initial shift $x_0^n$ is normalized by requiring that
\begin{equation}\label{trans-un00}
u_n(0,0)=\theta_0.
\end{equation}
\begin{lemma}\label{lem-trans-x0}
There exists $x_0^n$ so that $u_n(t,x)$ satisfies (\ref{trans-un00}) and
$\lim_{n\to+\infty}x_0^n=-\infty$. Moreover, there exists $N_0$ and $\eps>0$ so that
$x_0^n<-\eps n$ for all $n>N_0$.
\end{lemma}
{\bf Proof.} Let $v_n(t,x;y)$ be the solution of (\ref{eq0}) with the
initial data $v_n(t=-n,x;y)=\zeta(x-y)$ -- we are looking for $x_0$
such that $v_n(0,0;x_0)=\theta_0$.  Note that (\ref{trans-monotone})
implies $v_n(0,0;0)>\theta_0$. In addition, the function
$\psi(x)=\exp(-\lambda(x-ct))$ is a super-solution for (\ref{eq0})
provided that
\begin{equation}\label{trans-cl}
c\lambda\ge \lambda^2+Mg^{max},
\end{equation}
with the constant $M>0$ chosen so that $f_0(u)\le Mu$. Let us choose $\lambda>0$
and $c>0$ sufficiently large so that (\ref{trans-cl}) holds.
The maximum principle
implies that there exists a constant $C>0$ so that
\[
v_n(t,x;y)\le C\exp\{-\lambda (x-y-c(t+n))\},
\]
and thus
\[
v_n(0,0;y)\le C\exp\{\lambda (y+cn)\}\le \frac{\theta_0}{2}
\]
for $y<-cn-K$ with a constant $K>0$. By continuity of $v_n(0,0;y)$ as
a function of $y$ there exists $x_0\in(-cn-K,0)$ such that
$v_n(0,0;x_0)=\theta_0$. In order to see that $x_0^n\to-\infty$ as $n\to +\infty$,
observe that $v_n(t,x;y)\ge w_n(t,x;y)$, where $w_n(t,x;y)$ is the
solution of the Cauchy problem
\begin{equation}\label{trans-vmin}
\pdr{w_n}{t}=\frac{\partial^2w_n}{\partial x^2}+g^{min}f_0(w_n),~~w_n(-n,x;y)=\zeta(x-y).
\end{equation}
Note that if $y$ stays uniformly bounded from below as $n\to+\infty$:
$y\ge K$ for all $n$, then, as in (\ref{trans-toone}), $w_n(0,0;y)\to
1$, which contradicts (\ref{trans-un00}), thus $x_0^n\to-\infty$ as
$n\to+\infty$. The refined estimate $x_0^n<-\eps n$ follows the
results of \cite{Rjfr} on the exponential in time convergence of the
solution of (\ref{trans-vmin}) to a sum of two traveling waves of
(\ref{trans-vmin}) moving with a positive speed $c^{min}>0$ to the right
and left, respectively. In particular, this implies that if
$y>-nc^{min}/2$ then for $n$ sufficiently large we have
$w_n(0,0;y)>(1+\theta_0)/2>\theta_0$ which is a contradiction. $\Box$

The standard elliptic regularity estimates imply that the sequence of
functions $u_n(t,x)$ is uniformly bounded, together with its
derivatives and thus the limit $u(t,x)=\lim_{n\to+\infty}u_n(t,x)$ is
a global in time and space, monotonically increasing solution to
(\ref{eq0}). The main remaining difficulty is to show that $u(t,x)$
has the correct limits as $x\to\pm\infty$ and its ``interface width''
is uniformly bounded in time so that it is indeed a transition front
in the sense of Berestycki and Hamel. The rest of the proof of
Theorem~\ref{thm-transition} is based on the following estimates for
any solution of the Cauchy problem (\ref{eq0}) with the initial data
$u(0,x)=\zeta(x-x_0)$, with any $x_0\in\Rm$ (we set here the initial
time $t_0=0$ for convenience).

\subsection*{The interface width estimate}

For $h \in (\theta_0,1)$ and $k \in (0,\theta_0)$, let $X_h^l(t)$
and $X_k^r(t)$ be defined by
\br
&&X_h^l(t)=\max \left\{ x > x_0 \;|\;\; u(t,x') > h
\;\;\forall x' \in [x_0,x)\right\} \no \\
&&X_k^r(t) = \min \left\{ x > x_0 \;|\;\; u(t,x') < k
\;\;\forall x' \in (x,\infty)\right\}
\label{thetahkdef}
\er
Our goal is to show that the width of the front can be bounded by a
universal constant depending only on $g^{max}$ and $g^{min}$.
\begin{prop} \label{thm:khseparation}
Let $u(t,x)$ be a solution of (\ref{eq0}) with the initial data
$u(0,x)=\zeta(x-x_0)$ for some $x_0\in\Rm$. For any $h \in
(\theta_0,1)$ and $k \in (0,\theta_0]$, there are constants $K_h \ge
0$ and $C \geq 0$ depending only on $h$, $k$, $g^{max}$ and $g^{min}$ such that for
any $t>K_h$ we have $u(t,x_0)>h$, and
\begin{equation}\label{trans-intwidth}
0 < X_k^r(t) - X_h^l(t) \leq C<+\infty
\end{equation}
for all $t > K_h$. We can take $K_h=0$ for $h\in(\theta_0,h_0)$.
\end{prop}
Let us note that the time delay $K_h$ is introduced simply because
initially the solution may be below $h$ everywhere so that $X_h^l(t)$ is not defined
for small times.

\subsection*{The interface steepness bound}

The next crucial estimate provides a lower bound for the steepness of
the interface.  First, we use the following lemma to define the
interface location.
\begin{lemma}\label{lem:thetadef}
Let $u(t,x)$ be a solution of (\ref{eq0}) with the initial data
$u(0,x)=\zeta(x-x_0)$ for some $x_0\in\Rm$. For all $t > 0$, there
exists a continuous function (the right interface) $X(t)$, $t\ge 0$,
monotonically increasing in $t$ and satisfying: $x_0 < X(t)$,
$u(t,X(t)) = \theta_0$ and $u(x,t) < \theta_0$ for all $x >X(t)$,
$u(t,x) > \theta_0$ for all $x \in (x_0,X(t))$.
\end{lemma}
{\bf Proof:} This follows from the strict monotonicity of $u$ with
respect to time and the maximum principle which precludes $X(t)$ from
having jumps since $f(x,u)=0$ for $0\le u\le\theta_0$. $\Box$

\begin{prop}\label{thm:uxlower}
Let $u(t,x)$ be a solution of (\ref{eq0}) with the initial data
$u(0,x)=\zeta(x-x_0)$ for some $x_0\in\Rm$. Then the following hold.

(i) There is a universal constant $p > 0$ depending only on
$g^{max}$ and $g^{min}$ such that
\be
u_x(t,X(t)) < - p \label{uxlower1}
\ee
for all $t > 0$, and
\be
u(t,x + X(t)) \leq \theta_0 e^{- p x} \label{uxlower2}
\ee
for all $x > 0$, $t > 0$.

(ii) There exists a universal constant $\delta>0$ so that
\be
u_t(t,X(t)) > \delta \label{utlower}
\ee
for all $t > 1$.  Moreover, for any $t_1 > 0$,
there are constants $H > 0$ and $L>0$ such that
\begin{equation}\label{trans-interf-speed}
\hbox{$0<L < \dot X(t) < H<+\infty$ for all $t \geq t_1$.}
\end{equation}
The constants $L$ and $H$ also depend only on the constants
$g^{min}$, $g^{max}$ and the
function $f_0$.
\end{prop}

\subsection*{The end of the proof of Theorem~\ref{thm-transition}}

Theorem~\ref{thm-transition} is an immediate consequence of
Propositions~\ref{thm:khseparation} and~\ref{thm:uxlower}. Consider
the sequence of functions $u_n(t,x)$ defined for $t\ge -n$ as
solutions of the Cauchy problem for (\ref{eq0}) with the initial data
(\ref{trans-indata}) and $x_0^n$ fixed by normalization
(\ref{trans-un00}). As we have mentioned above, the standard elliptic
regularity estimates imply that there exists a subsequence
$n_k\to+\infty$ so that $u_{n_k}(t,x)$ converge locally uniformly,
together with its derivatives, to a limit $u(t,x)$ which is a global
in time and space solution to (\ref{eq0}), monotonically increasing in
time. Moreover, the interface locations $X_n(t)$ converge to $X(t)$
such that $u(t,X(t))=\theta_0$, $0<L<\dot X(t)<H$, and $X(0)=0$.  The
normalization (\ref{trans-un00}) implies that $u(0,0)=\theta_0$, and,
in addition, the bounds (\ref{uxlower1})-(\ref{utlower}) hold for the
limit $u(t,x)$. The upper bound (\ref{uxlower2}) implies immediately
that $u(t,x+X(t))\to 0$ as $x\to+\infty$ uniformly in $t$. It remains only to check
that $u(t,x+X(t))\to 1$ as $x\to-\infty$ uniformly in $t$. To this
end, assume that there exists $\eps_0>0$ and a sequence of points
$t_m\in\Rm$, and $x_m\to-\infty$ such that
$u(t_m,x_m+X(t_m))<1-\eps_0$. This, however, contradicts
(\ref{trans-intwidth}) with $k=\theta_0$ and $h=1-\eps_0$.  Therefore,
$u(t,x)$ is, indeed, a transition wave. $\Box$

\commentout{

\begin{lemma}\label{lem:vdef}
Let $u_n(t,x)$ be a solution of (\ref{eq0}) with the initial data
(\ref{trans-indata}) which satisfies the normalization
(\ref{trans-un00}), and $n \geq 1$. There exists a constant $p > 0$, an increasing function $\beta(t)$, and a nonincreasing in $x$ function $v(x)$ such that $v(0) =
\theta_0$, $v_x(0) <- p$, $\lim_{t \to \infty} \beta(t) = +\infty$,
\begin{eqnarray}
&&\lim_{x \to - \infty} v(x) =  1, \label{trans-vlower}\\
&&\lim_{x \to + \infty} v(x)  =  0,\label{trans-vupper}
\end{eqnarray}
and
\br
&&u_n(t,x + X_n(t)) \geq v(x), \quad \forall \;x \in [- \beta(t + n), 0],
\label{un-v1} \\
&&u_n(t,x + X_n(t))  \leq  v(x), \quad \forall \;x \in [0, \infty],
\label{un-v2}
\er
for all $t\ge 0$. The constant $p$ and the functions $v(x)$ and $\beta(t)$ depend only on $g^{max}$ and $g^{min}$.
\end{lemma}

}

A convenient way to restate some of the properties of the functions $u_n(t,x)$ we
will need later is as follows.  Let us define the class of admissible non-linearities
\[
{\cal G}=\{f(x,u)=g(x)f_0(u):~g^{min}\le g(x)\le g^{max},~g(x)\in C(\Rm)\}.
\]
\begin{lemma}\label{lem:vdef}
Given $0<g^{min}\le g^{max}<+\infty$ and $f_0(u)$, there exists $p > 0$
and a non-increasing in~$x$ function $v(x)$, such that
$v(0) =\theta_0$, ${v'(0)} <- p$,
\begin{eqnarray}
&&\lim_{x \to - \infty} v(x) =  1, \label{trans-vlower}\\
&&\lim_{x \to + \infty} v(x)  =  0,\label{trans-vupper}
\end{eqnarray}
and the following holds: given any solution $u_n(t,x)$ of (\ref{eq0}) with
$f(x,u)\in{\cal G}$, and with the initial data
(\ref{trans-indata}) which satisfies the normalization
(\ref{trans-un00}),
and for any $R > 0$, we have
\br
&&u_n(t,x + X_n(t)) \geq v(x), \quad \forall \;x \in [-R, 0],
\label{un-v1} \\
&&u_n(t,x + X_n(t))  \leq  v(x), \quad \forall \;x \in [0, \infty],
\label{un-v2}
\er
for all $t\ge 0$, if $n$ is sufficiently large, depending only on
$R$, $g^{min}$, and $g^{max}$. The function $v(x)$ depends only
on the constants $g^{max}$ and $g^{min}$, and the function $f_0(u)$.
\end{lemma}
{\bf Proof:} Setting $v(x)=\theta_0 e^{-px}$ for $x\ge 0$, with $p>0$
as in Proposition~\ref{thm:uxlower} we see that the upper bound
(\ref{un-v2}) follows from (\ref{uxlower2}), and (\ref{trans-vupper})
is obviously satisfied, as well as a strictly negative upper bound for
$v'(0)$ .

In order to define $v(x)$ for $x<0$ we consider a solution of (\ref{eq0})  with
$f\in{\cal G}$, which satisfies~(\ref{trans-un00}),
and with initial data as in (\ref{trans-indata}),
and choose $t_1=1$ and
find the corresponding $L$ as in Proposition~\ref{thm:uxlower}, so that $\dot X_n(t)\ge L$ for $t \ge -n+1$. Now, $X_n(t)\ge x_n^0+ L(t+n-1)$,
and thus for $t\ge 0$, $n \geq N_R=1 + R/L $, and
$x\in[-R,0]$ we have
\begin{eqnarray*}
&&x+X_n(t)\ge x_n^0 + L(t+n-1)- R \geq x_n^0.
\end{eqnarray*}
For $h \in [\theta_0,1)$, let $X^l_{n,h}(t)$ be defined by
(\ref{thetahkdef}). We use the convention that $X^l_{n,h}(t) =
-\infty$ if $u_n(t,x_n^0) < h$. For any $h' \in [\theta_0,1)$,
$X^l_{n,h}(t)$ is finite for all $t > 0$ and for all $h \in
[\theta_0,h']$, if $n>N(h')$ is sufficiently large, depending only on
$g^{min}$ and $h'$.  This follows directly from
Proposition~\ref{thm:khseparation}.
Now, for $x < 0$ and $n \geq 1$,
define
\be
v_n(x;f) = \sup \left\{ h \in [\theta_0,1): \;\; \sup_{t \ge 0}
\left(X_n(t) - X^l_{n,h}(t) \right) \leq \abs{x} \right\}. \label{vndef}
\ee
We indicate above explicitly the dependence of $v_n$ on the nonlinearity
$f(x,u)$.  Then we set
\[
v(x) = \inf_{f\in{\cal G}}\inf_{n\ge 1+|x|/L}  v_n(x;f).
\]
The set of possible values of $h$ over which the supremum is taken in
(\ref{vndef}) contains $\theta_0$. Therefore, $v(x) \geq \theta_0$ for
all $x<0$. From (\ref{vndef}) it is easy to see that $v_n(x;f)$ is
non-increasing in $x$ (for $x < 0$) for each $f\in{\cal G}$, and $v_n(0;f) = \theta_0$.
Hence $v(x)$ is also non-increasing in $x$ and $v(0) = \theta_0$.

Next we show that $v(x) \to 1$ as $x \to -\infty$. For any $h \in [\theta_0,1]$,
Proposition~\ref{thm:khseparation} implies that for all $f\in{\cal G}$ we have
\[
\sup_{t \ge 0} \left(X_n(t) - X^l_{n,h}(t) \right) < C(h)
\]
for some finite constant $C(h)$,  depending only on $g^{min}$ and
$g^{max}$, provided that $n>N(h)$, which ensures that $u_n(0,x_n^0)\ge h$.
So, for $x$ such that both $x < -C(h)$ and $1+|x|/L>N(h)$,
we have $v(x) > h$. Since $h$
may be chosen arbitrarily close to $1$, it follows that
$\lim_{x \to -\infty} v(x) = 1$.

Finally, in order to see that (\ref{un-v1}) holds, fix $R > 0$ and $f\in{\cal G}$,
and let $u_n(t,x)$ be the solution of the corresponding Cauchy problem.
By definition of $v$,
\[
\hbox{$v(x) \le v_n(x;f)$ for all  $x \in [-R,0]$,}
\]
provided that $n\ge 1+R/L$.
Therefore,
\be
X_n(t) - X^l_{n,h}(t)  \leq - x
\ee
for all $h \in [\theta_0,v(x)]$, $n\ge 1+R/L$, and all $t \geq 0$.
Hence, $X^l_{n,h}(t) \geq x + X_n(t) \geq x_0^n$ so that
\be
u_n(t,x + X_n(t)) \geq h \hbox{ for all $h \in [\theta_0,v(x)]$},
\ee
and for all $t \geq 0$ and $n\ge 1+R/L$. This proves (\ref{un-v1}). $\Box$



\subsection*{Bounds for the location of level sets}

In order to finish the proof of Theorem~\ref{thm-transition} it
remains to prove Propositions~\ref{thm:khseparation}
and~\ref{thm:uxlower}. We need first to establish some simple bounds
on the location of the level sets of the function $u(t,x)$. Let
$c^{min}$ and $c^{max}$ be the unique speeds of the traveling wave
solutions of the constant coefficient equations
\[
-cq_x=q_{xx}+f^{min}(q),~q(-\infty)=1,~q(+\infty)=0,
\]
and
\[
-cq_x=q_{xx}+f^{max}(q),~q(-\infty)=1,~q(+\infty)=0,
\]
respectively. The next lemma will allow us to relate the position $X_h^l(t)$ to
$X_{h'}^l(t-s)$ with $s>0$ and $h'<h$ -- this allows us to control the width
of the front in the back, where $u$ is close to $1$.
\begin{lemma} \label{lem:thetahlower}
Let $\delta > 0$ and let $u(t,x)$ be as in
Propositions~\ref{thm:khseparation} and~\ref{thm:uxlower}.  Suppose,
in addition, that $u(0,x) > \delta + \theta_0$ for all $x \in
[x_L,x_R]$.  If $\sigma = \abs{x_R - x_L}$ is sufficiently large,
depending on $\delta$, then for any $h \in (\theta_0, 1)$ there are
constants $\beta > 0$ and $\tau > 0$, such that
\be\label{trans-Xh-below}
X_h^l(t) \geq x_R + c^{min}t  - \beta
\ee
for all $t \geq \tau$. The constants $\beta$ and $\tau$ depend only on
$h$, $\delta$, $\sigma$, and $f^{min}$.
\end{lemma}
{\bf Proof:} This follows from the comparison principle and the
stability results of \cite{Rjfr}. Specifically, if $\sigma$ is
sufficiently large, consider the function $v(x,t)$ which solves the
equation
\[
v_t = \Delta v + f^{min}(v)
\]
with the initial data
\[
v(0,x)=(\delta+\theta_0)\chi_{[x_L,x_R]}(x)
\]
at $t = 0$. Then, as we have mentioned, the results of \cite{Rjfr}
imply that $v$ converges as $t\to+\infty$ to a pair of traveling waves
moving to the left and the right with speed $c^{min} > 0$. The
convergence is exponentially fast. Therefore, after some time $\tau$,
which depends on $h$ and on the convergence rate, $v(t,x) \geq h$ on
the set $[x_R,x_R + c^{min}t - \beta]$, for some constant $\beta >
0$. The maximum principle implies that $u(t,x) \geq v(t,x)$ and
(\ref{trans-Xh-below}) follows.  $\Box$

\begin{corollary} \label{cor:thetalower}
Let $h \in (\theta_0, 1)$. There is a constant $\beta$ depending only on
$h$ and $f^{min}$ such that for $t_1 \geq 0$ sufficiently large,
\be
X_h^l(t) \geq X_h^l(t_1) + \max \left( 0, c^{min} (t - t_1) - \beta \right)
\ee
for all $t \geq t_1$.
\end{corollary}
{\bf Proof:} Because $u$ is strictly monotone increasing in $t$,
$u(t,x_0) > h$ and $X_h^l(t)$ is well-defined for $t > t_1$, if $t_1$
is sufficiently large. After this time, $X_h^l(t)$ is an increasing
function of $t$. The bound now follows from Lemma~
\ref{lem:thetahlower} with $\delta = h-\theta_0$, $x_L = x_0$, $x_R =
X_h^l(t_1)$, and replacing $t=0$ with $t_1$. We may choose $t_1$
larger, if necessary, so that $\abs{x_R - x_L} = \abs{X_h^l(t_1) -
x_0} \geq \sigma$. $\Box$

\begin{lemma} \label{lem:thetaupper}
Suppose that $u(0,x) \leq C e^{- c^{max} x}$ for all $x$.
Then there is a constant $\eta > 0$ such that
\be
X(t) \leq c^{max}t  + \eta.
\ee
\end{lemma}
{\bf Proof:} This follows from the comparison principle and the
stability results of \cite{Rjfr}.  \qed

From these bounds it follows immediately that
\be
c^{min} \leq \liminf_{t \to \infty} \frac{X(t)}{t }
\leq \limsup_{t \to \infty} \frac{X(t)}{t}  \leq c^{max}. \label{cpmbound}
\ee

\subsection*{The proof of Proposition~\ref{thm:uxlower} (i)}

The strategy of this proof is a version of the sliding
method~\cite{BNsliding}. Suppose $\eps \in (0, c^{min})$ is
sufficiently small so that initially at time $t=0$ we have
\be
u_x(0,X(0)) \leq - \eps \theta_0 \label{thmux1}
\ee
and
\be
u(0,x + X(0)) \leq \theta_0 e^{- \eps x}, \label{thmux2}
\ee
both for $x\ge 0$.
The proof is in two steps. First,
we claim that there is $\tau_1 > 0$ such that
\br
&&u_x(t,X(t)) \leq  - \eps \theta_0 \label{thmux3} \\
&&u(t,x + X(t)) \leq  \theta_0 e^{- \eps x} \label{thmux4}
\er
hold for all $t \in [0,\tau_1]$ and $x\ge 0$. Then, replacing $t=0$ by
$\tau_1$ in (\ref{thmux1}) and (\ref{thmux2}), we may proceed
inductively, so that the bounds (\ref{thmux3}) and (\ref{thmux4}) also
hold for all $t \in [\tau_1, \tau_2]$, for some $\tau_2 >
\tau_1$. Continuing in this manner we obtain an increasing sequence
$\{ \tau_k \}_{k=1}^\infty$. In the second step we will show that the
increments $\tau_{k+1} - \tau_k$ must be uniformly bounded away from
zero. Therefore, $\sup_k \{ \tau_k \} = +\infty$, so that
(\ref{thmux3}) and (\ref{thmux4}) must hold for all $t \ge 0$.

Now we establish the claim.  Fix $y > 0$, and define
\be
\psi(t,x;y) = \theta_0 e^{-\eps(x - \eps t - y - X(0))}. \no
\ee
Then \[
\psi_t - \Delta_x \psi = 0
\]
and $\psi(0,X(0) + y) = \theta_0$.  Since $f(x,u) = 0$ for $x > X(t)$, the
function $w = \psi - u$ satisfies
\[
w_t - \Delta w = 0
\]
in the region
\be
\{ (t,x)\;|\;\; x \geq X(t), \; t \ge 0 \}. \no
\ee
Moreover, $w(0,x) > 0$ for all $x \geq X(0)$, and $w(t,X(t)) > 0$ if
$t>0$ is sufficiently small.  Let $\tau_y $ be the first time after
$t=0$ at which $w(t,x) = 0$ for some point $x \geq X(t)$:
\be
\tau_y = \sup \{ t \in (0, \infty) \;|\;\; w(t,x) > 0
\;\forall x \in [X(\tau),\infty),\;\forall \tau \in [0,t) \}.
\ee
Since $\eps < c^{min}$, the lower bound in (\ref{cpmbound}) and the
definition of $\psi$ implies that $\tau_y$ must be finite. The maximum
principle implies that the value $w(\tau_y,x) = 0$ is attained at the
boundary: $w(\tau_y,X(\tau_y)) = 0$. Therefore, at $t = \tau_y > 0$,
we have
\be
\psi(\tau_y,X(\tau_y)) = u(\tau_y,X(\tau_y)) = \theta_0,
\ee
and thus
\[
X(\tau_y) = \eps\tau_y + y + X(0).
\]
Using this identity in the inequality
$u(\tau_y,x + X(\tau_y)) \leq  \psi(\tau_y,x + X(\tau_y))$
for $x\ge 0$, we obtain
\be
u(\tau_y,x + X(\tau_y)) \leq  \psi(\tau_y,x + X(\tau_y)) =
\theta_0 e^{-\eps x}\hbox{ for $x\ge 0$.} \label{uuppertau}
\ee
Moreover, since $w(\tau_y,x) \geq 0$ for all $x \geq X(\tau_y)$, with
equality at $x = X(\tau_y)$, we have
\be
u_x(\tau_y,X(\tau_y)) \leq \psi_x(\tau_y,X(\tau_y)) =
- \eps \theta_0. \label{uxuppertau}
\ee
The time $\tau_y > 0$ is a continuous, increasing function of $y$ such
that $\lim_{y \to 0^+} \tau_y = 0$. Therefore, if we fix $y>0$ then
for any $t<\tau_y$ we may find $y'$ such that
$\tau_{y'}=t$. Therefore, (\ref{uuppertau}) and (\ref{uxuppertau})
together imply that
\begin{eqnarray*}
&&u(t,x + X(t))  \leq  \theta_0 e^{- \eps x},  \\
&&u_x(t,X(t)) \leq  - \eps \theta_0
\end{eqnarray*}
both hold for all $t \in [t_0,\tau_y]$. This establishes the first
claim with $\tau_1 = \tau_y$.

Why must $\tau_{k+1} - \tau_k$ be bounded away from zero, uniformly
for all $k \geq 1$? Since $u_x(t,X(t)) < - \eps \theta_0$ for all $t
\in [\tau_1, \tau_{k+1}]$, regularity estimates imply that
\be\label{trans-Xdot-up}
\dot X(t) = - \frac{u_t(t,X(t))}{u_x(t,X(t))} \leq
\frac{\norm{u_t(\cdot, t)}_\infty}{\eps \theta_0} \leq H
\ee
for some constant $H < \infty$ for all $t \in [\tau_1, \tau_{k+1}]$.
However, at $t = \tau_{k+1}$, we have
\be
\eps(\tau_{k+1} - \tau_k) + y + X(\tau_k) = X(\tau_{k+1})
\leq X(\tau_k) + H(\tau_{k+1} - \tau_k).
\ee
Therefore,
\be
(\tau_{k+1} - \tau_k) \geq \frac{y}{H - \eps} > 0.
\ee
Since we may fix $y$, independently of $k$, this lower bound is
independent of $k$. We have now proved (\ref{uxlower1}) and
(\ref{uxlower2}). The upper bound $\dot X(t) \leq H$ now follows from
(\ref{uxlower1}), and the regularity estimates for $u$, as in
(\ref{trans-Xdot-up}). This finishes the proof of part (i) of
Proposition~\ref{thm:uxlower} and of the upper bound for $\dot X(t)$
in part (ii) of this proposition.

\subsection*{The proof of Proposition~\ref{thm:khseparation}}

The standard elliptic regularity estimates imply that for any
$t_1 > 0$, there exists a constant $M>0$ so that
$\norm{u_{xx}(t,\cdot)}_\infty < M$ for all $t > t_1$. Therefore, we
have for $x < X(t)$, using (\ref{uxlower1}):
\br
&&u(t,x)\geq u(t,X(t)) + u_x(t,X(t))(x - X(t)) - \frac{1}{2} M (x - X(t))^2 \no \\
&& ~~~~~~~~\geq  \theta_0 - p (x - X(t)) - \frac{1}{2} M (x - X(t))^2 \label{quadlower}
\er
for all $t > t_1$. For $x = X(t) - p/M$, this gives,
\be
u(t,x) \geq \theta_0 + \frac{p^2}{2M}.\label{ulowerp2M}
\ee
Since $u$ is monotone in time, this implies that there is $t_2 \geq t_1$
such that for all $t \geq t_2$,
\be
X_h^l(t) \geq X(t) - p/M,  \label{thetaHbnd1}
\ee
if $h \in (\theta_0, \theta_0 + p^2/2M]$. The time gap between $t_1$
and $t_2$ may be needed to allow $u(t,x)$ go get above the value $h$
on the interval between the points $x_0$ and $X(t_1)-p/M$.

The exponential bound (\ref{uxlower2}) implies that for any $k \in (0,
\theta_0)$, the distance
\[
0 < X_k^r(t) - X(t) \leq  \frac{1}{p} \abs{\log k - \log \theta_0}
\]
is bounded uniformly in time. Combining this with (\ref{thetaHbnd1}),
we see that there is a constant $C_1$ such that for any $h \in
(\theta_0, \theta_0 + p^2/2M]$ and any $k\in(0,\theta_0)$, we have
\be
X_k^r(t) - X_h^l(t) \leq C_1
\ee
for all $t \geq t_2$.

Now suppose that $h \in (\theta_0+p^2/2M,1]$. First, Lemma~\ref{lem:monzeta}
implies that there exists a time $t_3$ so that $u(t_3,x_3)>h$ at some
point $x_3>x_0$ -- hence, the function $X_h^l(t)$ is well-defined for
$t\ge t_3$.  In addition, there exists a constant $C>0$ so that
$X(t_3)\le x_0+C$ at this time.  Let $\delta = p^2/2M$,
$\gamma=\theta_0+\delta$, and define $d_0 > 0$ by
\[
d_0 = \sup_{t \in [0, t_2 + \tau]} \left(X_\gamma^l(t) - X_h^l(t)\right)
\leq X(t_2 + \tau) - x_0 \leq C_2 + c^{max}(t_2 + \tau)<+\infty
\]
where $\tau$ is the constant from Lemma \ref{lem:thetahlower} with
$x_L = x_0$ and $x_R = X_\gamma^l(t_2)$. If necessary, we may take
$t_2$ to be larger so that $\abs{x_R - x_L} = \sigma$ is sufficiently
large according to Lemma \ref{lem:thetahlower}. Now for any $t \geq
t_2 + \tau$, we apply Lemma \ref{lem:thetahlower} with the starting
time $t_0 = t - \tau$, $\delta = p^2/2M$ and $x_R = X_\gamma^l(t -
\tau) \geq X_\gamma^l(t_2)$. We conclude that
\[
X_\gamma^l(t) - X_h^l(t)\leq X_\gamma^l(t) - X_\gamma^l(t - \tau) - c^{min}\tau
+ \beta.
\]
As we have already shown in the proof of Proposition
\ref{thm:uxlower}(i), we have $\dot X(t) \leq H$ for all $t >
t_1$. Therefore, since $t - \tau \geq t_2 \geq t_1$, we have
\br
&&X_\gamma^l(t) - X_h^l(t)  \leq  X_\gamma^l(t) - X_\gamma^l(t - \tau) -
c^{min}\tau + \beta
\leq X(t) - X_\gamma^l(t - \tau) - c^{min}\tau + \beta \no \\
&&=   \left(X(t) - X(t - \tau)\right) + \left(X(t - \tau) -
X_\gamma^l(t - \tau) \right) - c^{min}\tau + \beta \no \\
&& \leq
H \tau + \left(X(t - \tau) - X_\gamma^l(t - \tau) \right) - c^{min}\tau + \beta
\leq   H \tau + C_1 - c^{min}\tau + \beta. \no
\er
This holds for all $t \geq t_2 + \tau$ Therefore, for any $t \geq 0$, we have
\[
X_\gamma^l(t) - X_h^l(t) \leq d_0 + H \tau + C_1 - c^{min}\tau + \beta
\]
so, the conclusion of Propostion~\ref{thm:khseparation} holds. \qed

\subsection*{The proof of Proposition~\ref{thm:uxlower}(ii)}

We may now prove (\ref{utlower}), the lower bound on
$u_t(t,X(t))$. Fix $h \in (\theta_0,h_0)$, $k = \theta_0$, and let $C
= C(h,k)$ be the constant from Proposition~\ref{thm:khseparation}. For
any $t_1 < t$ (but sufficiently away from $t=0$), we have $\abs{\dot
X(s)} \leq H$ for all $s \in [t_1,t]$, so
\[
X(t_1) \in [ X(t) - H  (t- t_1), X(t)].
\]
It follows from Corollary \ref{cor:thetalower} that there is a
constant $\beta$, independent of $t_1$, such that
\be
X_h^l(t_1 + \Delta t) \geq X_h^l(t_1)
+ \max \left( 0 , c^{min}(\Delta t) - \beta \right)
\label{thetaHC2}
\ee
for all $\Delta t > 0$. So, if we choose, $\Delta t = (C + \beta)/c^{min}$,
(\ref{thetaHC2}) implies
\[
X_h^l(t_1 + \Delta t)  \geq X_h^l(t_1) + C \geq X(t_1).
\]
The last inequality follows from Proposition~\ref{thm:khseparation}
and our choice of $C$. Therefore, we have
\[
u(t_1 + \Delta t,X(t_1)) \geq h,
\]
and by the Mean Value Theorem there must a point $t_2 \in [t_1,t_1 +
\Delta t]$ such that
\[
u_t(t_2,X(t_1)) \geq (h - \theta_0)/(\Delta t),
\]
since $u(t_1,X(t_1)) = \theta_0$.

Now, let $t_1 = t - 2 \Delta t$ (recall that $\Delta t$ is defined
independently of $t_1$). Thus there exists a point $x_2 \in [ X(t) - 2
H \Delta t, X(t)]$ and $t_2 \in [t - 2\Delta t, t - \Delta t]$ such
that $u_t(t_2,x_2) \geq r > 0$, where $r = (h - \theta_0)/(\Delta t)$.

The function $q(t,x) = u_t(t,x)$ satisfies a PDE of the form
\[
q_t = \Delta q + V(x,t) q
\]
with $q \geq 0$ and $\norm{V}_\infty < \infty$. The Harnack inequality \cite{KrySaf}
implies that there is $K$ depending only on $H$, $\Delta t$, and
$\norm{V}_\infty$ such that
\[
q(t,X(t)) \geq K
\sup_{\substack{x_2 \in [ X(t) - 2 H \Delta t, X(t)]\\
t_2 \in [t - 2\Delta t, t - \Delta t] }} q(t_2,x_2) \geq K r > 0.
\]
Therefore, there is $\delta = K r$ depending only on the properties of
$f$ such that $u_t(t,X(t)) \geq \delta$ for all $t$ sufficiently
large. Since $u_t(t,X(t)) > 0$ for all $t > 0$, this implies
(\ref{utlower}).  Finally, the lower bound $\dot X(t) > L > 0$ now
follows from (\ref{uxlower1}), (\ref{utlower}), the first equality in
(\ref{trans-Xdot-up}), and the elliptic regularity estimates for
$u$. \qed


\section{Asymptotic spreading for the Cauchy problem} \label{sec:asympspread}

\subsection*{Spreading of monotonically increasing in time solutions}

Now we return to equation (\ref{eq0}) with a random reaction term, and
we prove Theorem~\ref{theo:asympspeed}. We first prove the result for
monotone increasing solutions. Consider the solution to (\ref{eq0})
with initial data $u_0(x,\omega) = \zeta(x + z_1)$ at time
$t=0$. Recall from the definition of the function $\zeta$ that
$\zeta(z_1) = \theta_0$ and $\zeta(x) < \theta_0$ for $x >
\theta$. Hence, we have $u_0(0,\omega) = \theta_0$. The initial data
looks like a bump-function with the right interface at the origin.  For
each realization $\omega \in \Omega$ of the random medium, the
following hold:
\bi
\item The solution $u(t,x,\omega)$ is strictly monotone increasing in $t$ and
all the estimates of Section~\ref{sec:tran} hold $\Pm$-a.s.
\item The function $X^+(t,\omega)$ defined by $u(t,X^+(t,\omega),\omega) = \theta_0$
and $X^+ \geq 0$ is well defined and continuous.
This defines uniquely the position of the
right-moving interface.
\item There are positive constants $C_{min}$, $C_{max}$, independent of
$\omega$ such that for $t > 1$ we have
$C_{min} \leq \dot X^+(t,\omega) \leq C_{max}$ .
\item For any $\xi \geq 0$, the time at which ``the interface reaches $\xi$",
denoted by $T(\xi,\omega)$, is well defined:
\be
\xi = X^+(T(\xi,\omega),\omega).
\ee
\ei
The first claim above follows from Lemma~\ref{lem:monzeta}, the second one
is a consequence of the maximum principle and monotonicity of $u(t,x,\omega)$
in time.
The last two claims are implied by (\ref{trans-interf-speed}).
Similarly, we may define the position $X^-(t,\omega)$ of the left-moving
interface by $u(t,X^-(t,\omega),\omega) = \theta_0$ and $X^-(t,\omega)\le -2z_1$
for $t\ge 0$.

The following proposition is a version of
Theorem~\ref{theo:asympspeed} for such monotonically increasing in
time solutions. We will then use a comparison argument to generalize
this result to arbitrary non-negative compactly supported initial
data as claimed in Theorem~\ref{theo:asympspeed}.
\begin{prop}\label{theo:cconv0}
There are nonrandom constants $c^*_+ \in [c^{min}, c^{max}]$ and
$c^*_- \in [-c^{max},-c^{min}]$ such that
\br
\lim_{t \to \infty} \frac{X^+(t,\omega)}{t} & = & c^*_+, \label{xplusconv}  \\
\lim_{t \to \infty} \frac{X^-(t,\omega)}{t} & = & c^*_- \label{xminusconv}
\er
hold almost surely with respect to $\Pm$, and in $L^1(\Omega,\Pm)$. For
any $\epsilon > 0$,
\be
\lim_{t \to \infty} \inf_{c \in [c^*_- + \epsilon,c^*_+ - \epsilon]}
u(t,ct,\omega) = 1 \label{uinnerlim}
\ee
and
\be
\lim_{t \to \infty} \sup_{c \in (-\infty,c^*_- - \epsilon]
\cup [c^*_+ + \epsilon,\infty)}  u(t,ct,\omega) = 0 \label{uouterlim}
\ee
hold almost surely with respect to $\Pm$, and in $L^1(\Omega,\Pm)$.
\end{prop}
\subsection*{Proof of Propositon~\ref{theo:cconv0}}

First, we explain that $X^+(t,\omega)$ is $\mathcal{F}$-measureable for each $t$.
Let $m$ be a positive integer. For each $m$ define the set of points
$ \{ x^m_j \} =  2^{-m} \mathbb{Z}$. For $m$ and $t$ fixed, let
\be
A^{m}_{j} = \{ \omega \in \Omega \;|\;\; u(t,x,\omega) \leq \theta_0, \;\forall x \geq x^m_j \}. \no
\ee
This is an $\mathcal{F}$-measureable set, since it is a closed set in
$C(\mathbb{R};[g^{min},g^{max}])$ (in the uniform convergence norm).
Define the random variable
\be
\eta^{m}(\omega) = \min_j \left( x_j^m \chi_{A^{m}_j}(\omega) \right)
\ee
where $\chi$ is the characteristic function. Since there are countably many
terms in the minimization, this is an $\mathcal{F}$-measureable random variable.
By definition, $ X^+(t,\omega) \leq \eta^m(\omega) \leq X^+(t,\omega) + 2^{-m}$.
Also, $\eta^{m}$ is nonincreasing in $m$. Therefore
\be
X^+(t,\omega)  = \lim_{m \to \infty} \eta^{m}(\omega)
\ee
and this must be $\mathcal{F}$-measureable, since the limit of a sequence of measurable functions is also measureable.

Next, we prove (\ref{xplusconv}) by using the sub-additive ergodic
theorem. Let us drop the superscript and denote $X(t,\omega) =
X^+(t,\omega)$. Given a positive integer $m\in\Nm$, let
$u^{(m)}(t,x,\omega)$ be the solution to (\ref{eq0}) for $t \ge 0$
with shifted initial data $u^{(m)}(x,0,\omega) = \zeta(x + z_1 - m)$ -- its
right interface is located initially at $x=m$. Let $X_m(t,\omega)\ge
m$, $t\ge 0$, denote the position of the corresponding right-moving
interface: $u^{(m)}(t,X_m(t,\omega),\omega) = \theta_0$. By
Proposition~\ref{thm:uxlower}, $X_m(t,\omega)$ satisfies the same
properties as $X(t,\omega)$, listed above. For $\xi \geq m$, let
$T_m(\xi,\omega)\ge 0$ denote the inverse of $X_m(t,\omega)$:
$u^{(m)}(T_m(\xi,\omega),\xi)=\theta_0$.

Now, for a pair of non-negative integers $m,n\in\Nm$, $n\ge m$,
define the family or random variables
\[
q_{m,n}(\omega) = T_m(n,\omega)
\]
which is the first time the interface hits the position $n$, when
started from position $m$. It is easy to see that for any integer $h
\geq 1$, the following translation invariance holds:
\be
q_{m + h,n + h}(\omega) = q_{m,n}(\pi_h \omega). \label{qstatprop}
\ee
The key observation in the proof of Proposition~\ref{theo:cconv0}
is the following ``near-subadditivity" lemma.
\begin{lemma} \label{lem:almostsubadd}
There exists a constant $\alpha > 0$ independent of $\omega$ such that
\be
q_{m,r}(\omega) \leq q_{m,n}(\omega) + q_{n,r}(\omega) + \alpha \label{alphasubadd}
\ee
holds for all pairs of integers $0 \leq m < n < r$.
\end{lemma}
We postpone the proof of this lemma for the moment and proceed with
the proof of Proposition~\ref{theo:cconv0}.  Using
Lemma~\ref{lem:almostsubadd} we now show that there is a nonrandom
constant $\bar q$ such that the limit
\[
\lim_{n \to \infty} \frac{1}{n} q_{0,n}(\omega) = \bar q
\]
holds almost surely.  Lemma \ref{lem:almostsubadd} shows that the
family $\{ q_{n,m} \}$ is ``almost" subadditive.  In order to turn it
into a truly sub-additive family define a new family
\[
\hat q_{m,n} = q_{m,n} + \beta (n - m)^{1/2}
\]
with $\beta$ sufficiently large to be chosen. The point here is that
$\hat q_{m,n}$ is a sub-linear correction of $q_{m,n}$. It also
preserves translation invariance of $q_{m,n}$: for any integer $h >
0$, we have, using (\ref{qstatprop}):
\[
\hat q_{m+h,n+h}(\omega)= q_{m+h,n+h}(\omega)  + \beta (n - m)^{1/2}
= q_{m,n}(\pi_h \omega)  + \beta (n - m)^{1/2} = \hat q_{m,n}(\pi_h \omega).
\]
Let $\alpha>0$ be as in (\ref{alphasubadd}) and choose $\beta>4\alpha$.
Then for any integers $0 \leq m < n < r$
the following elementary inequality holds:
\[
\alpha + \beta (r - m)^{1/2} - \beta (r - n)^{1/2}  - \beta (n - m)^{1/2} \leq 0
\]
since $r-n \geq 1$ and $n-m \geq 1$.  Lemma~\ref{lem:almostsubadd}
implies that with this choice of $\beta$ the family $\hat q_{m,n}$ is
sub-additive: for any integers $0 \leq m < n < r$ we have
\begin{eqnarray*}
&&\hat q_{m,r} = q_{m,r} + \beta (r - m)^{1/2}
\leq  q_{m,n } + q_{n,r} + \alpha + \beta (r - m)^{1/2}  \\
&&~~~~~= \hat q_{m,n } + \hat q_{n,r} + \left(\alpha +
\beta (r - m)^{1/2} - \beta (r - n)^{1/2}  - \beta (n - m)^{1/2}\right)
\le  \hat q_{m,n } + \hat q_{n,r}.
\end{eqnarray*}

Corollary \ref{cor:thetalower} implies that $\hat q_{m,r}$ is at most
linear: $0 \leq \hat q_{m,r} \leq C(1 + (m-r))$ for some constant
$C>0$. As the group $\pi_n$ acts ergodically on $\Omega$, we can apply
the subadditive ergodic theorem (see, e.g. \cite{Lig}) to conclude
that
\be
\lim_{n \to \infty} \frac{1}{n} \hat q_{0,n} =
\inf_{n > 0} \frac{1}{n} \Em\left[\hat  q_{0,n}\right] = \bar q
\ee
holds almost surely, where $\bar q$ is a deterministic constant.
By definition of $\hat q$, this implies that
\[
\lim_{n \to \infty} \frac{1}{n}  q_{0,n} = \bar q
\]
also holds almost surely. Since $q_{0,n} = T(n,\omega)$ and $X(t)$ is
increasing in $t$, it is easy to see that, as a consequence,
\[
\lim_{t \to \infty} \frac{X(t,\omega)}{t} = (\bar q)^{-1}:=  c^*_+
\]
holds almost surely. The fact that $c^*_+ \in [c^{min}, c^{max}]$
follows from (\ref{cpmbound}). This proves (\ref{xplusconv}), and the
proof of (\ref{xminusconv}) is identical.

The fact that limits (\ref{uinnerlim}) and (\ref{uouterlim}) hold is
an immediate consequence of (\ref{xplusconv}), (\ref{xminusconv}) and
the fact that the width of the interface is bounded by a universal
constant, as stated in Proposition~\ref{thm:khseparation}. This
completes the proof of Proposition~\ref{theo:cconv0}. \qed

\subsection*{The proof of Theorem~\ref{theo:asympspeed}}

Now, we use comparison arguments to extend
Proposition~\ref{theo:cconv0} to the case of any non-negative
deterministic initial data with a sufficiently large compact
support. By ``sufficiently large", we mean large enough so that the
solution does not converge uniformly to zero (extinction). Lemma
\ref{lem:monzeta} implies that the condition $u_0(x) \geq
\zeta(x-x_0)$ with some $x_0\in\Rm$ is sufficient to guarantee that
extinction does not occur.

Let $w_0(x)$ be compactly supported with $0 \leq w_0 \leq 1$ and
deterministic. Suppose that
\[
w_0(x) \geq \zeta(x-x_0)
\]
for some $x_0 \in \mathbb{R}$ and let $w(t,x,\omega)$ solve
(\ref{eq0}) with initial data $w_0(x)$. For each $t > 0$, let
$X^+(t,\omega)$ be the largest real number satisfying
$w(t,X^+(t,\omega),\omega) = \theta_0$.

If $u(t,x,\omega) $ solves the equation with initial data
$u(0,x,\omega)=\zeta(x-x_0) \leq w_0(x)$, Proposition~\ref{theo:cconv0}
applies to $u(t,x,\omega)$, and the maximum principle implies that
$w(t,x,\omega) \geq u(t,x,\omega)$. Therefore, $w(t,x,\omega)$
satisfies
\[
\lim_{t \to \infty} \inf_{c \in [c^*_- + \epsilon,c^*_+ - \epsilon]}
w(t,ct,\omega) \geq 1.
\]
Since $w \leq 1$ for all $t\ge 0 $, this implies the first bound of
Theorem \ref{theo:asympspeed}.

For the other bound, observe that for every realization $\omega$ we
have $\max_{x\in\Rm} w(t=1,x,\omega)<c_0<1$ with a deterministic
constant $c_0$.
The estimates in the previous section imply that there is a finite
time $\tau > 0$ depending only on the properties of $f$ such that
\[
w(t=1,x,\omega) \leq u(t=1+\tau,x,\omega) , \quad \forall x \in \mathbb{R}.
\]
Then the maximum principle implies that $w(s,x,\omega) \leq u(s + \tau,x,\omega)$
for all $s \geq 1$. This implies
\be
\lim_{t \to \infty} \sup_{c \in (-\infty,c^*_- - \epsilon]
\cup [c^*_+ + \epsilon,\infty)}  w(t,ct,\omega) \leq 0. \no
\ee
Since $w \geq 0$ for all $t$, this completes the proof of Theorem
\ref{theo:asympspeed}.
\qed

\subsection*{The proof of Lemma~\ref{lem:almostsubadd}}

Translation invariance (\ref{qstatprop}) implies that it suffices to
prove that (\ref{alphasubadd}) holds for $m=0$. We first show that
there is an integer $K > 0$ independent of $\omega$ such that for all
$r,n$ a ``delayed" version
\be
q_{0,r}(\omega) \leq q_{0,n}(\omega) + q_{n - j,r}(\omega) \label{njclaim}
\ee
holds for $j = \min(K,n)$.  Let $h = \max_x \zeta(x) \in (\theta_0,
1)$ and define $X_h^l(t)$ as in (\ref{thetahkdef}). By
Proposition~\ref{thm:khseparation}, there is a constant $C > 0$,
independent of $\omega$ such that
\begin{equation}\label{alm-xhx}
X_h^l(t) \geq X(t,\omega) - C.
\end{equation}

Now let $K$ be the smallest integer greater than $C + z_2+z_1$ (recall
that $\zeta(x) = 0$ for all $\abs{x} \geq z_2$). First,
(\ref{njclaim}) obviously holds for $n\le K$ as for such $n$ it
becomes
\[
q_{0,r}(\omega) \leq q_{0,n}(\omega) + q_{n - n,r}(\omega) =
q_{0,n}(\omega) + q_{0,r}(\omega),
\]
which is true since $q_{0,n}(\omega)\ge 0$.

If $n \geq K$ then (\ref{alm-xhx}) implies that
\[
u(T(n,\omega),x,\omega) \geq h, \quad \forall \;x
\in (-z_1,n-C)\subseteq
(-z_1,n - K + z_2).
\]
On the other hand, we have
\[
\hbox{$\zeta(x+z_1-(n-K))=0$ for $x\notin(-z_1,n-K+z_2)$.}
\]
Therefore, if $n \geq K$, we have
\[
u(T(n,\omega),x,\omega) \geq \zeta(x + z_1 - (n - K)) = u^{(n-K)}(0,x,\omega),
\hbox{ for all $x\in\Rm$.}
\]
Since the equation is invariant with respect to $t$, the maximum
principle implies that for any $s \geq 0$,
\[
u(T(n,\omega) + s,x,\omega) \geq u^{(n-K)}(s,x,\omega),
\]
thus $X(T(n,\omega)+s,\omega) \geq X_{n-K}(s)$. Now setting
$s = T_{n-j}(r,\omega) = T_{n-K}(r,\omega)$ we see that
\[
X(T(n,\omega) + T_{n-K}(r,\omega)) \geq X_{n-K}(T_{n-K}(r,\omega)) = r.
\]
Since $X$ is increasing in $t$, this implies $T(r,\omega) \leq
T(n,\omega) + T_{n-K}(r,\omega)$ which establishes (\ref{njclaim}) for
$n \geq K$.  Thus, the claim holds for all $n > 0$.

Using the fact that $u^{(n-j)}$ is
monotone in $t$ and the estimates of the previous section, one can
show that there is a constant $\alpha > 0$ independent of $n$ and
$\omega$ such that
\[
u^{(n-j)}(t,x,\omega) \geq \zeta(x + z_1 - n), \quad \forall x \in \mathbb{R}, \;
t \geq \alpha,
\]
where $j = \min(K,n)$ is bounded independent of $n$ and $\omega$.
This and the maximum principle imply that
\[
u^{(n-j)}(\alpha + s,x,\omega) \geq u^{(n)}(s,x,\omega),
\quad \forall x \in \mathbb{R}, \; s \geq 0.
\]
Thus, we have
\[
q_{n - j,r}(\omega) \leq q_{n,r}(\omega) + \alpha
\]
This inequality and (\ref{njclaim}) imply the desired result:
\[
q_{m,r}(\omega) \leq q_{m,n}(\omega) + q_{n,r}(\omega) + \alpha .
\]
This finishes the proof of Lemma~\ref{lem:almostsubadd}. \qed

\section{Random Traveling Waves} \label{sec:rtw}

Now we use the results of the previous sections to construct a random
traveling wave solution to the equation (\ref{eq0}) and prove Theorem
\ref{theo:rtwexist} and Corollary \ref{cor:profinv}.

\subsection{Construction of the Traveling Wave}

The starting point comes from the proof of Theorem A(1) in
\cite{Shen1}.  We consider a family $\tilde u_n(t,x,\omega)$ of
solutions of the Cauchy problem (\ref{eq0}) with the initial data
$\tilde u_n(t=-n,x,\omega)=\zeta^s(x-\tilde x_0^n(\omega))$. Here
$\zeta^s$ is the step function:
\[
\zeta^s(x) =  \left \{ \begin{array}{ll} 1, & x < 0, \\
0, & x \geq 0, \end{array} \right.
\]
and the shift $\tilde x_0^n(\omega)$ is fixed by the normalization, as
in (\ref{trans-un00})
\[
\tilde u_n(0,0,\omega)=\theta_0,~~\tilde u_n(0,x,\omega)<\theta_0\hbox{ for $x>0$}.
\]
In this section we denote with tilde objects related to solutions with step-like initial
data, while those without tilde correspond to those arising from bum-like
initial data.

The random initial shift $\tilde x_0^n(\omega)$ is measureable with respect to
$\mathcal{F}$ and is uniquely defined. The existence and uniqueness of
$\tilde x_0^n(\omega)$  follows from the fact that if $y_1 < y_2$, the
comparison principle implies that the solution to (\ref{eq0}) with initial data
$\zeta^s(x - y_1)$ must be below the solution with initial data $\zeta^s(x - y_2)$. Therefore, for fixed $n$, the front position at time $t=0$ is a monotonic function
of the shift, and the maximum principle implies that it is continuous. Then,
using arguments similar to those in the proof of Lemma \ref{lem-trans-x0} one
can show that there must be a unique $\tilde x_0^n(\omega) \in [-cn,cn]$ such
that the normalization condition is satisfied, if $c > 0$ is sufficiently large.

The measureability of $\tilde u_n$ and $\tilde x_0^n$ may be proved as
in~\cite{Shen1} (Theorem A(1), therein). For the readers' convenience we sketch
the proof now. For each $n$, let $w(t,x,\omega;y)$ solve (\ref{eq0}) for $t > -n$
with initial data $w(t = -n,x,\omega) = \zeta^s(x - y)$. Let $\eta_n(y,\omega)$ denote
the largest real number satisfying $w(0,\eta_n,\omega) = \theta_0$.  For each $y$,
$\eta_n(y,\omega)$ is $\mathcal{F}$-measureable. This may be proved as in the
case of $X^+(t,\omega)$ in Section \ref{sec:asympspread}. Now we vary $y$, and
we wish to choose $y$ as a measurable function of $\omega$ so that
$\eta_n(y,\omega) = 0$. For each positive integer $k$ define
$\{y^k_l\} = 2^{-k} \mathbb{Z}$.  Let $r$ be a positive integer, and define
\be
A^{k,r}_{l} = \{ \omega \in \Omega \;|\;\; \abs{\eta(y^k_l,\omega)} \leq  1/r\}. \no
\ee
This is an $\mathcal{F}$-measureable set since $\eta(y,\cdot)$ is
$\mathcal{F}$-measureable. Then we set
\br
\hat x^n_0(\omega) & = & \lim_{r \to \infty} \lim_{k \to \infty} \min_{l} \left( y^k_l \chi_{A^{k,r}_l}(\omega)  \right). \label{hatx}
\er
Notice that
\[
\min_{l} \left( y^k_l \chi_{A^{k,r}_l}(\omega)  \right)
\]
is $\mathcal{F}$-measureable, being the infimum of a countable set of
measurable functions, and it is nonincreasing in $k$ and nondecreasing in $r$.
Thus, the limits in (\ref{hatx}) exist and $\hat x^n_0(\omega)$ is measurable.
The continuity of $\eta(y,\omega)$ with respect to $y$ and the uniqueness of
$\tilde x_0^n$ imply that $\hat x_0^n(\omega) = \tilde x_0^n(\omega)$.
So, $\tilde x_0^n$ is $\mathcal{F}$-measureable.

The measureability of $\tilde u_n$ now follows from the measureability of
$\tilde x_0^n$. Specifically, for fixed $n$ and $t$, the function $\tilde u_n$
may be expressed as a composition of measureable maps:
\be
\tilde u_n(t,\,\cdot \,,\omega) = G_2 \circ G_1(\omega)
\ee
where $G_1(\omega):(\Omega,{\cal F}) \to
(\mathbb{R} \times \Omega,{\cal B}\times{\cal F})$ is the measureable
map $G_1(\omega) = (\tilde x^n_0(\omega),\omega)$ and
$G_2(y,\omega):(\mathbb{R} \times \Omega,{\cal B}\times{\cal F})\to
C(\mathbb{R};[0,1])$ is the measureable map defined by solution of (\ref{eq0}) with
initial data $\zeta^s(x - y)$ (shifted by $y$) at time $t = -n$. Here ${\cal B}$ is the Borel
$\sigma$-algebra on $\Rm$.

Now, for $\tilde
x_0^n(\omega)$ defined in this way, we wish to take a limit $n\to
+\infty$ to construct a global-in-time solution. That is, we wish to
define
\be
\tilde w(t,x,\omega) = \lim_{n\to +\infty} \tilde u_n(t,x,\omega), \label{Udef}
\ee
and show that this is a traveling wave solution. The existence of a
measureable limit, converging locally uniformly, and satisfying the
PDE follows from Shen \cite{Shen1} (see proof of Theorem A(1)) and
regularity estimates. A key observation in \cite{Shen1}, is that the
convergence (\ref{Udef}) holds as $n\to+\infty$, not just along a
particular subsequence $n_k$. This is because the functions $\tilde
u_n$ satisfy the following monotonicity relation at $t=0$:
\begin{eqnarray}
&&\tilde u_n(0,x,\omega) > \tilde u_m(0,x,\omega), \;\;\text{if}\;\;  x < 0\nonumber\\
 &&\tilde u_n(0,x,\omega)< \tilde u_m(0,x,\omega),\;\;\text{if}\;\;  x > 0 ,
 \label{wtildemonotone}
\end{eqnarray}
almost surely, for any $m >n$. Therefore, the function $\tilde
w(t,x,\omega)$ is measureable in $\omega$. However, the difficulty is
that the limit might be trivial: one may obtain $\tilde w(t,x,\omega)
\equiv \theta_0$ for all $x$ and $t$. Here is where we invoke the
results of the previous sections.

\subsection*{Uniform limits at infinity}

Using Proposition~\ref{thm:khseparation} and the estimates of
Section~\ref{sec:tran}, we can show that the limit $\tilde w$ must be
non-trivial.
\begin{lemma}\label{lem-tw-limits}
Let $\tilde w(t,x,\omega)$ be constructed as above.
Then  we have
\br
&&\lim_{x \to \infty} \sup_{\omega \in \Omega} \tilde w(t=0,x,\omega)  =  0, \no \\
&&\lim_{x \to -\infty} \inf_{\omega \in \Omega} \tilde w(t=0,x,\omega) =  1.
\label{wtildelimits}
\er
\end{lemma}
{\bf Proof.} We prove (\ref{wtildelimits}) by comparing the functions
$\tilde u_n(t,x,\omega)$ with functions $u_n(t,x,\omega)$ defined as
follows.  For each $n$, let $u_n(t,x,\omega)$ denote the solution of
(\ref{eq0}) with initial data $\zeta(x - x_0^n)$ at time $t = -n - 1$
(note that $u_n$ starts at time $t=-n-1$, and not at $t=-n$).  The
function $\zeta(x)$ is the bump-like sub-solution used in
Section~\ref{sec:tran}, so the solution $u_n(t,x,\omega)$ is strictly
monotone increasing in $t$ and the estimates of Section \ref{sec:tran}
apply to $u_n$.  The point $x_0^n= x_0^n(\omega)$ is a random shift
depending on $n$. For such initial data, let $X_n^+(t;x_0^n,\omega)$
be defined as in Lemma~\ref{lem:thetadef}. The random shift
$x_0^n(\omega)$ is chosen so that $X_n^+(0;x_0^n,\omega) = 0$ for
all $n\in\Nm$, $\omega \in \Omega$. This is the same normalization as
applied to $\tilde u_n(t,x,\omega)$. Existence of the shift
$x_0^n(\omega)$ for each realization $\omega$ follows from
Lemma~\ref{lem-trans-x0}.


Having defined the function $x_0^n(\omega)$, one can show that
for each $t > -n$, there exists a unique point $\xi_n(t,\omega)$ such
that
\begin{eqnarray*}
&&\tilde u_n(t,x,\omega) > u_n(t,x,\omega),\;\;\text{if}\;\;
x < \xi_n(t,\omega) \\
&&\tilde u_n(t,x,\omega) <  u_n(t,x,\omega) , \;\;\text{if}\;\;  x > \xi_n(t,\omega).
\end{eqnarray*}
That is, the graphs of the two solutions $\tilde u_n$ and $u_n$ intersect
at time $t$ only at the point $x = \xi_n(t,\omega)$. This may be
proved as in Lemma~4.6 of \cite{Shen1} using the results of
Angenent~\cite{Ang} and the maximum principle. Here we sketch the argument. Recall that despite the
suggestive notation we have initialized $u_n(t,x,\omega)$ at time
$t=-n - 1$ so that at time $t = -n$, we have $0<u_n(t=-n,x,\omega)<1$
everywhere. Therefore, using the approximation argument employed in
the proof of Lemma 4.6 of \cite{Shen1}, one may argue as if the graphs
of $u_n(t=-n,x,\omega)$ and $\tilde u_n(t=-n,x,\omega) = \zeta^s(x -
\tilde x_0^n)$ intersect at only one point. Since the function $q = \tilde u_n - u_n$ satisfies a PDE of the form
\be
q_t = \Delta q + V(t,x)q \no
\ee
with $\norm{V}_\infty < \infty $, Theorems A and B of \cite{Ang} show that the
zero set of the function $q(t,x)$ is discrete and cannot increase. Therefore,
the graphs of $\tilde u_n$ and $u_n$ have only one intersection point for all $t > -n$.
We have chosen $x_0^n$ and
$\tilde x_0^n$ so that at $t = 0$, the graphs intersect at $x=0$:
$\tilde u_n(0,0,\omega) = \theta_0 = u_n(0,0,\omega)$ almost
surely. Therefore, $\xi_n(0,\omega) = 0$, and both
\[
\tilde u_n(0,x,\omega) > u_n(0,x,\omega), \quad x < 0
\]
and
\[
\tilde u_n(0,x,\omega) < u_n(0,x,\omega), \quad x > 0
\]
must hold, $\Pm$-a.s. for all $n\in\Nm$.

Passing to the limit $n\to
+\infty$, we see that for $x < 0$ we have a lower bound for $\tilde
w(0,x,\omega)$:
\[
\tilde w(0,x,\omega) \geq \liminf_{n \to +\infty} u_n(0,x,\omega):=v^-(x,\omega).
\]
It follows that from Lemma \ref{lem:vdef} that $v^-(x,\omega)$ has a
deterministic lower bound
\[
\lim_{x \to -\infty} v^-(x,\omega) \geq \lim_{x \to -\infty} v(x) = 1,
\]
which holds for all realizations $\omega$. Similarly, for $x > 0$, we
have an upper bound for $\tilde w(0,x,\omega)$:
\[
\tilde w(0,x,\omega) \leq \limsup_{n \to +\infty} u_n(0,x,\omega):= v^+(x,\omega)
\]
and, once again, by Lemma \ref{lem:vdef}, $v^+(x,\omega)$ has a deterministic
upper bound:
\[
\lim_{x \to +\infty} v^+(x,\omega) \leq \lim_{x \to +\infty} v(x) =  0,
\]
that holds for all $\omega$.  This proves that (\ref{wtildelimits})
holds uniformly in $\omega$. $\Box$

\subsection*{The translation property}

We have know shown that $\tilde w(t,x,\omega)$ satisfies properties (i)-(iv) in the
definition of a random traveling wave. Since the limit $\tilde w(t,x,\omega)$ is nontrivial, the position of the interface
$\tilde X(t,\omega)$ may be defined at
time $t$:
\be
\tilde X(t,\omega) = \max \{ x \in \mathbb{R} \;|\;\;
\tilde w(x,t,\omega) = \theta_0 \} . \label{Xchar}
\ee
The measureability of $\tilde X(t,\omega)$ may be proved as in the case of $\tilde x^n_0(\omega)$.

Finally we show that the translation property (v) holds. The argument here is similar to that in \cite{Shen1}; we sketch details for the readers' convenience. Notice that we have not needed to assume that the index $n$ is an integer. In fact, we may assume $n \in [1,\infty)$.  The key observation that leads to property (v) is that for any $m \geq 0$,
\be
\tilde u_n(m,x + \theta_n(m,\omega),\omega) = \tilde u_{n+m}(0,x ,\pi_{\theta_n(m,\omega) }\omega) \label{twpropeq}
\ee
must hold. Here, $\theta_n(m,\omega)$ is the position of the interface at time $t=m$, when the solution is initialized at time $t = -n$ (with initial data $\zeta^s(x - \tilde x^n_0$)).   One may think of $\pi_{\theta_n(m,\omega) }\omega$ as the ``current environment" associated with the ``current location" of the interface (i.e. $\theta_n(m,\omega)$) at time $t = m$. If at time $t=m$ the interface is at $x = \theta_n(m,\omega)$, then in the coordinate system shifted by $\theta_n(m,\omega)$ the interface is at the origin. So if we simply shift $x$ by $ \theta_n(m,\omega)$ and $t$ by $m$, equality (\ref{twpropeq}) follows from the definition of $\tilde u_n$ and $\tilde u_{n+m}$, the fact that $f(x + \theta_n(m,\omega),u,\omega) = f(x,u, \pi_{\theta_n(m,\omega) }\omega)$, and the fact that $\tilde x^n_0$ and $\tilde x^{n+m}_0$ are uniquely defined. In particular, the function
\be
v(t,x,\pi_{\theta_n(m,\omega) } \omega) := \tilde u_n(t + m,x + \theta_n(m,\omega),\omega)
\ee
satisfies the shifted equation
\be
v_t = \Delta v + f(x + \theta_n(m,\omega),v,\omega) = \Delta v + f(x,v,\pi_{\theta_n(m,\omega) }\omega)
\ee
with initial data $v(t = -n - m, x,\pi_{\theta_n(m,\omega) } \omega) = \zeta^s(x - \tilde x^n_0(\omega) + \theta_n(m,\omega) )$. Since $\tilde x^{n+m}_0(\pi_{\theta_n(m,\omega) } \omega)$ is uniquely defined, this is the same initial value problem solved by $\tilde u_{n+m}(t,x,\pi_{\theta_n(m,\omega) } \omega)$. Therefore, uniqueness implies $v = \tilde u_{n+m}$. So, (\ref{twpropeq}) holds.

%

By definition of $\tilde w$ and $\tilde X$, $\theta_n(m,\omega) \to \tilde X(m,\omega)$ as $n \to \infty$, and the left hand side of (\ref{twpropeq}) converges to
\be
\lim_{n \to \infty} \tilde u_n(m,x + \theta_n(m,\omega),\omega) = \tilde w(m,x + \tilde X(m,\omega),\omega). \label{claimunm}
\ee
We claim that as $n \to \infty$ the right hand side of (\ref{twpropeq}) converges to $\tilde w(0,x,\pi_{\tilde X(m,\omega)}\omega)$. To see this, we express the right hand side of (\ref{twpropeq}) in the reference frame corresponding to $\tilde X(m,\omega)$. Let $\omega_m = \pi_{\tilde X(m,\omega)} \omega$ and define
\be
z_{n+m}(t,x,\omega_m) = \tilde u_{n + m}(t,x +  \tilde X(t,\omega) - \theta_{n}(m,\omega), \pi_{\theta_n(m,\omega) }\omega). \no
\ee
Then $z_{n+m}$ satisfies
\be
z_t = \Delta z + f(x + \tilde X(t,\omega) - \theta_{n}(m,\omega), z, \pi_{\theta_n(m,\omega) }\omega) = \Delta z + f(x, z, \omega_m) \label{zeqn}
\ee
with initial condition $z_{n+m}(t = - n -m,x,\omega_m) = \zeta^s(z - z^n_0)$ where $z^n_0 = \tilde X(t,\omega) - \theta_{n}(m,\omega) - \tilde x^{n - m}_0$. However, the function $\tilde u_{n + m}(t,x,\omega_m)$
satisfies the same equation (\ref{zeqn}) with initial condition $\tilde u_{n + m}(t = -n -m, x,\omega_m) = \zeta^s(z - \tilde x^n_0(\omega_m))$. In general, $z^n_0 \neq \tilde x^n_0(\omega_m)$, but the maximum principle still implies that at time $t= 0$ either $z_{n+m}(0,x,\omega_m) > \tilde u_{n+m}(0,x,\omega_m)$ for all $x$, or $z_{n+m}(0,x,\omega_m) < \tilde u_{n+m}(0,x,\omega_m)$ for all $x$.  However, at time $t = 0$, $\tilde u_{n+m}(0,0,\omega_m) = \theta_0$, and $z_{n+m}(0,\theta_{n}(m,\omega) - \tilde X(t,\omega)  ,\omega_m) = \theta_0$. Since $\lim_{n \to \infty} \abs{ \theta_{n}(m,\omega) - \tilde X(t,\omega)} = 0$, one can use the maximum principle to show that in the limit $n \to \infty$, the two functions coincide:
\be
\lim_{n \to \infty} z_{n+m}(t,x,\omega_m) = \lim_{n \to \infty} \tilde u_{n + m}(t,x,\omega_m) \no
\ee
for all $x$ and $t$, as in Lemma 4.5(2) of \cite{Shen1}, since they both converge to $\theta_0$ at the point $x = 0$, $t=0$. By definition of $\tilde w$, the right hand side at $t=0$ is simply
\be
 \lim_{n \to \infty} \tilde u_{n + m}(0,x,\omega_m) = \tilde w(0,x, \pi_{\tilde X(m,\omega)} \omega). \no
\ee
This proves the claim (\ref{claimunm}) and establishes the translation property
\be
\tilde w(0,x, \pi_{\tilde X(m,\omega)} \omega) = \tilde w(m,x + \tilde X(m,\omega),\omega). \no
\ee
This completes the construction of the traveling wave.

For later use, let us note that the preceding proof shows that the
function $W(x,\omega) = \tilde w(0,x,\omega)$ satisfies
\br
&&W(x,\omega)  \geq  v(x), \quad \forall \;x < 0\no \\
&&W(x,\omega) \leq v(x), \quad \forall \;x > 0\no
\er
where $v(x)$ is deterministic and defined in Lemma
\ref{lem:vdef}. Therefore, the translation property (v) implies
that
\br
&&\tilde w(t,x + \tilde X(t),\omega)  \geq  v(x), \quad \forall \;x < 0\no \\
&&\tilde w(t,x + \tilde X(t),\omega) \leq  v(x), \quad \forall \;x > 0\no
\er
also holds.

\subsection*{Traveling waves and generalized transition waves}

Let us point out that an alternative way to establish existence of a
traveling wave is to use the bump functions $u_n(t,x,\omega)$ and pass
to the limit along a subsequence $n_k(\omega)\to+\infty$ to obtain a
non-trivial transition front $u(t,x,\omega)$ in the sense of
Berestycki and Hamel. Theorem A of \cite{Shen1} shows that a traveling
wave will exist if there exists such a generalized transition front
for each realization. However, it may be necessary to take the limit
along a different subsequence $n_k(\omega)$ for each $\omega$. This
may result in a transition wave $u(t,x,\omega)$ that may not be
measureable. The advantage of using a shift of the step function
$\zeta^s(x)$ is that the sequence is monotone in the sense of
(\ref{wtildemonotone}) and the limit (\ref{Udef}) may be taken as
$n\to+\infty$.  Therefore, the limit is measureable.

\commentout{
\begin{definition}[see also Shen \cite{Shen1}, Def. 2.3] \label{def:wls}A solution $\tilde v(x,t,\omega):\mathbb{R} \times [0,\infty) \times \Omega \to \mathbb{R}$ of (\ref{eq0}) is called a {\bf wave-like solution} if for any $h,k \in (0,1)$ with $h > k$,
\be
0 \leq \theta_{k}^+(t,\omega) - \theta_h^-(t,\omega) \leq C \label{wlsineq}
\ee
for all $t \geq 0$, where
\br
\theta_h^-(t,\omega) = \sup \{ x \in \mathbb{R} |\; \tilde v(x',t,\omega) > h \;\;\forall x' < x \} \no \\
\theta_{k}^+(t,\omega)=  \inf \{ x \in \mathbb{R} |\; \tilde v(x',t,\omega) < k \;\;\forall x' > x \}
\er
and $C = C(h,k)$ is a constant independent of $t$ and $\omega$.
\end{definition}
Roughly speaking, a wave-like solution is a solution for which there are uniform, global-in-time bounds on the width of the interface. The random traveling wave $\tilde w$ we have constructed is itself a wave-like solution, however, we were not able to construct a separate wave-like solution to aid in proving that the limit $\tilde w$ is non-trivial.

}

\subsection{Properties of the traveling wave} \label{sec:rtwproperties}

Now, we finish the proof of Theorem~\ref{theo:rtwexist} -- it remains
to show that the interface location $\tilde X(t)$ is a strictly
increasing function and that the limit in (\ref{rtwsamec}) exists and
is deterministic. First, we show that
\be
\lim_{t \to \infty} \frac{\tilde X(t,\omega)}{t} = c^*_+ \label{wtildeconv}
\ee
almost surely with respect to $\Pm$, where $c^*_+$ is the
deterministic right spreading rate defined in
Theorem~\ref{theo:asympspeed}.  Using Theorem~\ref{theo:asympspeed}
and the comparison principle, it is easy to show that
\[
\liminf_{t \to \infty } \frac{\tilde X(t,\omega)}{t} \geq
\liminf_{t \to \infty } \frac{ X(t,\omega)}{t} = c^*_+,
\]
with probability one, since we may construct compactly supported
initial data that fits below each realization of the profile
$W(x,\omega)$.

\subsection*{A super-solution for the traveling wave}

For an upper bound, we construct a super-solution related to a
construction in \cite{FMcL}. Let $u_n(t,x,\omega)$ be the same family
of monotone increasing solutions constructed in the proof of
Theorem~\ref{theo:rtwexist}.  Let $q \in (0,\theta_0/3)$ and set $h =
1 - q$. For $v(x)$ defined as in Lemma~\ref{lem:vdef}, let $y_h =
v^{-1}(h) < 0$ (i.e. $v(y_h) = h$). Pick $n \in \Nm$ sufficiently large so that
Lemma~\ref{lem:vdef} holds with $R = -y_h$. Therefore, by Lemma \ref{lem:vdef}, we have
\br
&&u_n(t,x + X_n(t),\omega) \geq  v(x), \quad \forall x \in [y_h, 0], \no \\
&&u_n(t,x + X_n(t),\omega) \leq  v(x), \quad \forall x > 0, \label{vboundu}
\er
for all $t \geq 0$. For a function $\gamma(t)$ to be chosen, define
\br
\bar u_n(t,x,\omega) = \left \{\begin{array}{ll} \min(1,u_n(\gamma(t),x,\omega) + q)
&  x > X_n(\gamma(t)) - y_h ,\\
1 & x \leq X_n(\gamma(t)) - y_h.   \end{array} \right.
\er
The function $\gamma(t)$ will be chosen so that $\gamma(0) > 0$ and
$\gamma'(t) > 1$.  We want to pick $\gamma(t)$ so that $\bar u$ is a
super-solution for $t \geq 0$. By construction, $\bar u$ now has a
wave-like profile, and $\bar u = 1$ for $x$ sufficiently negative.

If $u_n(\gamma(t),x,\omega) \geq h$ or $x < X_n(\gamma(t)) - y_h$, then
$\bar u_n(t,x,\omega) = 1 \geq \tilde w(t,x,\omega)$. On the other hand,
if $u_n(\gamma(t),x,\omega) \leq h$ and $x \geq X_n(\gamma(t)) - y_h$, then
$\bar u(t,x,\omega) \leq 1$ and
\br
\pdr{\bar u_n}{t} - \frac{\partial^2\bar u_n}{\partial x^2} - f(x,\bar u_n) =
\left( \gamma'(t) - 1 \right) \pdr{u_n}{t} + \left[ f(x,u_n) - f(x,\bar u_n) \right].
\label{Nuexp}
\er
Now we show that the right hand side of (\ref{Nuexp}) can be made
non-negative for $x \geq X_n(\gamma(t)) - y_h$, so that $\bar u_n$ is
a super-solution in this region.

By the properties of $f$, there exists $s \in (0,\theta_0/3)$ such
that $f(x,u) - f(x,\bar u) \geq 0$ wherever $ u \geq 1 - s$ and $\bar
u \leq 1$. Note that such an $s$ may be chosen independently of $q$
and $h$. For such an $s $ fixed, (\ref{vboundu}) and the properties of
$v$ imply that there is $\beta > 0$ such that
\be
\left \{ x \in [y_h,\infty)\; |\;\; u_n(\gamma(t),x + X_n(\gamma(t)),\omega)
\in [s, 1 - s] \right \} \subset [-\beta, \beta] \label{Rbound}
\ee
for all $t \geq 0$. By Proposition~\ref{thm:uxlower}, there is $\delta > 0$ such that
\[
\pdr{u_{n}}{t}(\gamma(t),X_n(\gamma(t)),\omega) > \delta.
\]
This and the Harnack inequality imply that there is $\eps > 0$ such
that
\begin{equation}\label{tw-harnack-ut}
\pdr{u_n}{t}(\gamma(t),x + X_n(\gamma(t)),\omega) > \eps,
\quad \forall x \in [-\beta,\beta], \;\;
t \geq 0.
\end{equation}

Now, if $x \in [X_n(\gamma(t)) - y_h,X_n(\gamma(t)) - \beta]$, then by
(\ref{Rbound}) we have $\bar u_n(t,x) \geq u_n(\gamma(t),x) \geq 1 -
s$, so $f(x,\bar u_n)\le f(x,u_n)$, the last term on the right side of
(\ref{Nuexp}) is non-negative and thus (\ref{Nuexp}) implies that in
this interval
\begin{equation}\label{tw-posit-rhs}
\pdr{\bar u_n}t -\frac{\partial^2 \bar u_n}{\partial x^2} - f(x,\bar u_n) =
\left( \gamma'(t) - 1 \right) \pdr{u_n}{t} \geq 0,
\end{equation}
since $\gamma'(t)\ge 1$.

If $x \in [X_n(\gamma(t)) + \beta,+\infty)$, then $u_n(\gamma(t),x) \leq
s$, so $\bar u_n(t,x) \leq s + q < \theta_0$. Hence $f(x,u_n) =
f(x,\bar u_n) = 0$ in this region, so again (\ref{tw-posit-rhs})
holds.

Finally, if $x \in [X_n(\gamma(t)) - \beta, X_n(\gamma(t)) + \beta]$, the
right side of (\ref{Nuexp}) can be bounded below using
(\ref{tw-harnack-ut}) by
\[
\partial_t \bar u - \bar u_{xx} - f(x,\bar u) \geq
\left( \gamma'(t) - 1 \right) \eps + \left[ f(x,u) - f(x,\bar u) \right]
 \geq \left( \gamma'(t) - 1 \right) \epsilon - K q \] where $K > 0$ is
 the Lipschitz constant for $f$. So if we choose $\gamma'(t) = 1 +
 {Kq}/{\eps}$, the right side is non-negative. For $\gamma(t)$ chosen
 in this way, we see that $\bar u_n$ is a super-solution wherever
 $\bar u < 1$, for all $t \geq 0$. Since $u_n$ is monotone increasing
 in $t$, we may also choose $\gamma(0)$ sufficiently large so that
\[
\bar u_n(0,x,\omega) \geq \tilde w(0,x,\omega).
\]
Therefore, the maximum principle implies that $\bar u_n(t,x,\omega)
\geq \tilde w(x,t,\omega)$ for all $t \geq 0$. Hence, we have
\[
\limsup_{t \to \infty} \frac{\tilde X(t)}{t} \leq
\limsup_{t \to \infty} \frac{ X_n(\gamma(t))}{t}=
\limsup_{t \to \infty}  \frac{ X_n(\gamma(t))}{\gamma(t)} \frac{\gamma(t)}{t}
=  c_+^* \left( 1 + \frac{Kq}{\eps}\right) =
c_+^* \left( 1 + \frac{K(1 - h)}{\eps}\right)
\]
Since $h$ can be chosen to be arbitrarily close to $1$, the right side
can be made arbitrarily close to $c_+^*$. Note that $s$ and $\beta$ can be chosen
independently of $h$, so that the parameter $\eps$ does not become
small as $h \uparrow 1$. This proves the upper bound and establishes
(\ref{wtildeconv}).

\subsection*{Monotonicity of the right interface}

We now prove the last claim of Theorem~\ref{theo:rtwexist} -- that the
interface $\tilde X(t)$ always moves to the right.
\begin{lemma}
For almost every $\omega \in \Omega$, the function $\tilde X(t,\omega)$ is
differentiable and strictly increasing in $t$.
\end{lemma}
\nit {\bf Proof:}
The maximum principle and the fact that $f(x,u) = 0$ for $u \leq \theta_0$
implies that $\tilde X$ cannot have jumps to the right:
\be
\limsup_{h \to 0^+} \tilde X(t + h,\omega) \leq \tilde X(t,\omega). \label{nojumps}
\ee
To see that $\tilde X$ is continuous and differentiable, note that
\br
\theta_0 = \tilde w(t,\tilde X(t),\omega) \label{wxdiffquotient}
\er
for all $t$. The function $W(x,\omega) = \tilde w(x,0,\omega)$ satisfies
\br
&&W(x,\omega)> v(x), \;\;\text{if}\;\;  x < 0\no \\
&&W(x,\omega)< v(x),  \;\;\text{if}\;\;  x > 0,\no
\er
$\Pm$-almost surely, and $v_x(0) < - p$ for some constant $p > 0$.
Therefore, we have
\[
W_x(0,\omega) = \tilde w_x(t,\tilde X(t),\omega) < - p < 0.
\]

The Implicit Function Theorem applied to (\ref{wxdiffquotient})
implies that there is a $C^1$ function $Y(t)$ such that $\theta_0 =
\tilde w(Y(t + h),t + h,\omega)$ for $h$ sufficiently small, and $Y(t)
= \tilde X(t)$. This, combined with the definition (\ref{Xchar}) and
(\ref{nojumps}), implies that $\tilde X(t)$ is continuous and that we
may differentiate (\ref{wxdiffquotient}) to obtain
\[
\tilde X'(t,\omega) = -
\frac{\tilde w_t(t,\tilde X(t,\omega),\omega)}
{\tilde w_x(t,\tilde X(t,\omega),\omega)}
 < \infty.
\]
This may also be written as
\br
\tilde X'(t,\omega) & = &  - \frac{W_{xx}(0, \pi_{\tilde X(t,\omega)} \omega) +
f(0,W(0, \pi_{\tilde X(t,\omega)} \omega ), \pi_{\tilde X(t,\omega)} \omega)}
{W_{x}(0, \pi_{\tilde X(t,\omega)} \omega)}. \no
\er

We have already shown that there is a set of full measure $\Omega_0
\subset \Omega$ such that $\Pm(\Omega_0) = 1$, and $\tilde
X(t,\omega)/t \to c^*_+ \geq c^{min} > 0$ for all $\omega \in
\Omega_0$ as $t\to+\infty$.  If $\tilde X(t)$ is not
strictly increasing in time, there are $t_1, t_2 \in \mathbb{R}$ such
$t_2 > t_1$ and $\tilde X(t_1,\omega_0) = \tilde X(t_2,\omega_0)$ for
some $\omega_0 \in \Omega_0$. Then
\[
\tilde w(t_1,x,\omega_0) =  W(x - \tilde X(t_1,\omega_0),
\pi_{\tilde X(t_1,\omega_0)} \omega_0) =
W(x - \tilde X(t_2,\omega_0), \pi_{\tilde X(t_2,\omega_0)} \omega_0)
= \tilde w(t_2,x,\omega_0)
\]
holds for all $x \in \mathbb{R}$. Hence, the function
$\tilde w(t,x,\omega_0)$ is periodic in $t$. This contradicts the fact that
$\tilde X(t,\omega)/t \to c^*_+ > 0$ for all $\omega \in
\Omega_0$. Therefore, $\tilde X(t + h,\omega) > \tilde X(t,\omega)$
for all $t \in \mathbb{R}$, $h > 0$, $\omega \in \Omega_0$. The proof
of Theorem~\ref{theo:rtwexist} is now complete. \qed

\subsection*{Proof of Corollary \ref{cor:profinv}}

This follows immediately from the definition of $\tilde X$ and $\tilde T$:
\[
\tilde w(\tilde T(\xi,\omega), x + \xi,\omega) =
W(x + \xi - \tilde X(\tilde T(\xi,\omega),\omega),
\pi_{ \tilde X(\tilde T(\xi,\omega),\omega)} \omega)
 =  W(x + \xi - \xi,\pi_{ \xi} \omega)
 =  W(x,\pi_{\xi}\omega).
\]
The last term on the right side is stationary with respect to shifts
in $\xi$ since the action of $\pi$ is measure-preserving. \qed

\bibliographystyle{plain}

\end{document}